\input amssym.def
\input amssym     
\magnification=1200
\baselineskip=16pt

\centerline{\bf Some New Applications of Weyl's Multi-Polarization Operators}
\bigskip
\centerline{\bf  Jacob Towber}  
\bigskip
\centerline{\it Department of Mathematics, De Paul University, Chicago,
 Illinois 60614}
\centerline{E-mail: jtowber@condor.depaul.edu}
\bigskip \bigskip
\centerline{\bf Introduction: Part One}
\bigskip
Let ${\frak A}_N$ denote the enveloping algebra of ${\frak gl}_N({\Bbb C})$.
\bigskip
In Chapter Two of {\it The Classical Groups},[W], Weyl constructed 
certain distinguished elements of ${\frak A}_N$, which he called
{\it quasi-compositions}; for a more modern presentation of this material,
see (Fulton and Harris,[FH,pp 514-515].). One 
of the main purposes of the present paper, is 
to direct the attention of representation-theorists to these remarkable
 (but rather neglected) objects in ${\frak A}_N$.

These objects will here be referred to, equivalently, as 
 ``{\it multi-polarizations}''
or as ``{\it Weyl polarizations}''; their definitions will be explained
 in \S1.3 below.
  Weyl used them for the purpose of proving
the celebrated Capelli identities in  ${\frak A}_N$, which in turn he 
utilized to prove what is often called (following Weyl's
nomenclature) the `First Fundamental Theorem of Invariant Theory' for the
classical groups. 

But, in fact, these Weyl polarizations are not necessary, (and are no longer 
 usually used), for either of these two original
purposes. Indeed, Capelli himself did
not require (nor mention) these multi-polarizations, at least not in
[Cap1]( the reference cited by Weyl in [Weyl], p.39), nor in 
[Cap2] (as cited by Young in [Young], pp.64--71), nor in  
[Cap3]. (For interesting summaries of Capelli's
work, cf. ([Young,loc.cit.], [Umeda].)  

Let it be emphasized at this point, that the
present paper is {\bf NOT} further concerned, either with the 
Capelli identities, or with the First Fundamental 
Theorem of Invariant Theory. 
\bigskip

Thus, to make plausible the claim (here proposed) that the Weyl polarizations
deserve further investigation, would seem to require obtaining some newer 
applications of these objects, to topics of more current interest, and with
results which have not been obtained by other methods.
Two such applications will be presented in the following paper, one to the
Verma-Shapovalov element (mapping one Verma module into another), the 
other to the remarkable complexes (constructed by A.Zelevinsky in [Zel]),
 which (together
 with their `conjugates', under the
interchange of symmetrization and alternation), furnish higher syzygies,
for the Weyl modules defined in [CL], and for the Pl\"ucker equations
which define the dual Weyl modules (also called shape functors or 
Schur functors).    

We postpone until Part Two of the introduction, a
sketch of issues involved in the second of these new applications of 
the Weyl polarizations. 

Here are the main results concerning the Weyl polarizations, to 
be established in the present paper:

I) Weyl's original construction of the Weyl polarizations, as  
differential operators, is recalled in Def.1.3.1 below. A 
remarkably simple and quite
different-looking, purely combinatorial construction (which, it is
hoped, is new) is given in Def.1.3.2; these two constructions are proved
equivalent in Th.2.3.1.
It is this combinatorial construction which plays a natural role in 
the study of the syzygies of the Pl\"ucker equations, as will be explained
 in \S3.1, and which make it plausible (it is here proposed) that still further
applications may exist for the Weyl polarizations.

II) The definition of the Verma-Shapovalov element is reviewed in
\S3.3 below. The work of Zelevinsky ([Zel]) and its later development by
  Akin ([Akin1,2]) relates this 
element to the above-mentioned complexes, and hence to the Weyl polarizations.
(I should like to thank Bhama Srinivasan for calling this connection
to my attention.) There is obtained below an explicit formula, stated in
\S3.3, expressing the action on
$$\underbrace{SV\otimes SV\otimes\cdots\otimes SV}_{N\;times}$$
of the Verma-Shapovalov element
for ${\frak sl}_N({\Bbb C})$, as a ${\Bbb Z}$-linear combination of
Weyl polarizations. (This appears not to be an obvious consequence of
two explicit formulas for the Verma-Shapovalov element which 
have appeared in the literature, 
but which do not utilize the Weyl polarizations,
 viz. the formula given by Carter and Lusztig
in [CL], and that given by[MFF].)

III) Explicit formulas in terms of the Weyl polarizations, which 
the author believes are new, are given
in \S3.1, for the differentials in the Zelevinsky complex discussed in 
the following sub-section.

These formulas are shown in \S3.2 to check with earlier precise data  
for the special case $N=3$, presented by Doty in ([Doty],p.134--136), and
 there attributed to Verma.

IV) The proofs of the results just cited, utilize formulas
(analogous to the Pieri formulas) for the product of a Weyl polarization
 by an elementary polarization.These formulas are derived in \S2.1
and \S2.2.
\medskip

The present author wishes to acknowledge his debt to D.-N. Verma, for 
explaining
in a number of public lectures over the past 20 years or so,  
the problem of explicitly defining the differentials
in the
above-mentioned complex; and to acknowledge in addition, a  number
of very illuminating conversations with Verma concerning these matters.
The author is also indebted to
J.Humphreys and A.Zelevinsky, for helpful suggestions concerning the
rather large literature involving the Shapovalov element.
 
Two members of the De Paul University Information Services who were especially
helpful in the preparation of this manuscript were Dan Wanek and Rick
Cruz. Thanks are also due the University of Illinois at Chicago for
extending the author its hospitality, as a Visiting Scholar, during the period
this paper was written. Gratitude is especially owed the UIC library
staff for help in obtaining some of the nineteenth-century work
of Capelli; I should especially like to thank D.L.Thomas (Library Associate
for the UIC Mathematics Department), and also Helen Badawi (Public Service 
Supervisor for the Science Library) and Gladys Odegaard (Reference Librarian
for  the Science Library).

Last (but not least), I wish to acknowledge my wife Diane's help
and encouragement.
\bigskip

\centerline{\bf Introduction: Part Two The Zelevinsky
Complex}
\medskip
The Weyl polarizations seem extremely well-suited, to furnishing
explicit formulas for the differentials in certain 
remarkable complexes which occur in representation-theory.
\medskip
The starting point here is the mid-nineteenth-century
Jacobi-Trudi identity between two symmetric polynomials in
$x_1,\cdots,x_m$, namely
$$s_\alpha =\hbox{ det }(h_{\alpha_i-i+j})=
\left|\matrix{h_{a_1}&h_{a_1+1}&\cdots&h_{a_1+N-1} \cr
             h_{a_2-1} & h_{a_2} & \cdots & h_{a_2+N-2} \cr
             \vdots&\vdots&\ddots&\vdots \cr
             h_{a_N-N+1}&h_{a_N-N+2}&\cdots&h_{a_N} \cr }
\right| \eqno(1)$$
Here $$\alpha = (a_1 \ge a_2 \ge \cdots \ge a_N \ge 0)$$
denotes a partition (possibly terminated by a string of zeros), $s_{\alpha}$
denotes the associated Schur function, and $h_i=h_i(x_1,\cdots,x_m)$ has its
usual meaning (i.e. the sum of all monomials in $x_1,\cdots,x_m$ of total
degree $i$; in particular, this is 1 if $i=0$, and 0 if $i<0$~).
\medskip
Let us denote the Grothendieck ring of polynomial
representation-functors of the group-scheme $GL_{\Bbb C}$  by
${GL_{\Bbb C}}^{\wedge}$ (Caution: The results in the present
paper relevant to the Zelevinsky complex, 
 are only asserted in characteristic 0---this is signalled by our
choice of ${\Bbb C}$ as groundfield. On the other hand, many of
the results obtained below concerning the Verma-Shapovalov
element, are valid for arbitrary ground-ring. Finally, only the
Lie algebra $A_N$ is here studied.)

A fruitful heuristic principle in the representation-theory
of $GL_{\Bbb C}$ and the
symmetric groups, is that whereby identities involving
symmetric polynomials can in many cases
be re-interpreted as equations in
the Grothendieck ring   ${GL_{\Bbb C}}^{\wedge}$.
Let us
now examine the result of applying this principle
to the equation (1) above:

In this heuristic manner, there is strongly suggested the
existance of a complex in the Grothendieck ring ${GL_{\Bbb C}}^{\wedge}$, which
is to be an exact 
sequence, such that equating to 0 its Euler-Poincare
characteristic, yields precisely the Jacobi-Trudi identity (1).
Over the last 25 years, a number of experts (discussed in further detail in 
\S3.1)
have considered the construction of such a complex. 
This project  seems to have been initiated by A. Lascoux in [Las];
 however, the complex constructed by Lascoux, unlike the 
 later Zelevinsky complex, has differentials which are basis-dependent,
(utilizing the combinatorics of Young tableaux) 
and so are not $GL(V)$-linear. These various attempts to construct
such a complex, have produced the following mixed results:

The {\it terms} in  these complexes may be precisely specified, in a simple
fashion for which all authors seem in agreement; the precise definitions of the
{\it differentials} in these complexes has been more problematical,
and it is (in the present author's opinion) only the work
of Zelevinsky in [Zel] 
which has first succeeded in producing mathematically acceptable
constructions of these differentials
as natural transformations. One of the main purposes of the
present paper, is that of rendering more transparent the construction
of these differentials, utilizing properties (to be established 
in Chapter Two below)
of Weyl polarizations. This work logically divides into two parts:

\noindent 1) We shall construct {\it explicitly} linear transformations
which play the role of the differentials in these complexes. This will
be done in \S3.1 below, and this construction uses nothing of the 
machinery of Verma modules and the homomorphisms between these.
These results appear to be quite new (except for the 
very special explicit results of Verma
 for the case $N=3$ of 3-part partitions, explained in Doty
[Doty,pp.134--136], which agree
with the present more general results---as will be verified below in
\S3.2 ).

\noindent 2)In Chapter Four, the explicit  maps furnished by this
'elementary' construction will be
proved to coincide precisely, with the maps provided
by Zelevinsky in terms of the theories
of Verma[Verma1,2], Bernstein-Gel'fand-Gel'fand[BGG]
 and Shapovalov[Shap]. To the 
present author's best knowledge,
it is only with the work  of Zelevinsky 
that such a comparison becomes possible---
there was in this sense no {\it specific} complex in ${GL_{\Bbb C}}^{\wedge}$ 
in existence earlier, for
which such comparison could make sense. It is for this reason that 
it seems appropriate to use 
here the terminology `Zelevinsky
complex'.  Some important later work
by Akin, Doty and Maliakas, developing these ideas of Zelevinky, may 
be found in [Akin1,2], [Doty2] and [Mal]. In particular, the present
paper is heavily indebted to the observation of Akin (in [Akin2,p. 418])
relating the differentials in the Zelevinsky complex to the Verma-Shapovalov
elements---this will be discussed further in \S3.3 below.
\bigskip
In the present author's opinion, there should exist a more 
direct combinatorial proof
that the complex constructed in \S3.1 is exact. At present, the 
only proof in the author's possession,
involves the detour (presented in Chapter Four below)
 through the theory of Verma modules
and the work of Zelevinsky.
One benefit of this detour is the explicit formula, presented in \S3.3
below, expressing (as mentioned earlier) the 
action on $S^{\otimes N}$ of each Verma-Shapovalov element for $A_N$ 
as a ${\Bbb Z}$-linear combination of Weyl polarizations.

\bigskip
\centerline{\bf Chapter 1  Weyl's Multi-Polarization Operators}
\medskip

\centerline{\S 1.1 \bf Notation}
\medskip
Throughout the following paper, the ground-field will be ${\Bbb C}$.
$\bigotimes$ will always denote $\bigotimes_{\Bbb C}$.
Let N be a positive integer.We denote by $\underline N$ the set
$\{ 1,...,N \}$
and by $S^{\otimes N}$ the functor (on the category of complex vector-spaces 
and ${\Bbb C}$-linear transformations) which assigns to each complex 
vector-space $V$, the complex vector-space 
$$S^{\otimes N}V=\underbrace{SV\otimes SV \otimes\cdots\otimes SV}_{\rm N \;
factors} $$ 
(where $SV$ is the usual symmetric algebra on $V$). Similarly, we define
$\Lambda^{\otimes N}$ to be the functor given by
$${\Lambda}^{\otimes N}V=\underbrace{{\Lambda}V\otimes {\Lambda} V 
\otimes\cdots\otimes {\Lambda}V}_{\rm N \;
factors} $$  
For $i,j \in \underline N$, we denote as usual by $E_{i,j}$ the $N \times N$ 
matrix whose only non-zero entry is a $1$ in row $i$
 and column $j$.
We denote by ${\frak A}_N$ the enveloping algebra of
  $\frak {gl}(N,{\Bbb C})$, and by $({\frak A}_N)^0$ the enveloping
algebra of $\frak {sl}(N,{\Bbb C})$ 
 (considered as a subalgebra of  ${\frak A}_N$).
In all computations in the present paper, permutations will act
on the {\it left} of elements, and linear transformations will act 
on the {\it left} of (column) vectors; thus compositions of
linear transformations and of permutations are to be read right-to-left.
\bigskip
\centerline{\bf \S 1.2 Elementary Polarization Operators}
\medskip
If $i$ and $j$ are integers between $1$ and $N$, there is
defined a natural transformation  
$$D_{i,j}:S ^{\otimes N}V \to S^{\otimes N}V \eqno (1.2.1)$$
(the {\it elementary polarization operator}), as follows.
(We review this well-known concept in some detail, in preparation for
its generalization in the next \S.)

We recall two definitions for these $D_{i,j}$, the first only meaningful
if V is finite-dimensional over a ground-field of 
characteristic $0$, the second valid for an arbitrary module V over an 
arbitrary commutative 
ground-ring --- and with both definitions equivalent for 
finite-dimensional V  over a field of characteristic 0.

Suppose first that V is a finite-dimensional complex vector-space, with 
${\Bbb C}$-basis 
$$ {\cal B}=(x_1,\ldots,x_M)$$ 
${\cal B}$ induces ${\Bbb C}$-algebra isomorphisms 
$$ \phi_{\cal B}:SV\simeq {\Bbb C}[X_1,\ldots,X_M] \eqno(1.2.2)$$
and
$$ \Phi_{\cal B}:S^{\otimes N}V\simeq {\Bbb C}[X_1^{(1)},\ldots,X_M^{(1)};
\ldots;X_1^{(N)},\ldots,X_M^{(N)}] \eqno(1.2.2a)$$
(where, in (1.2.2a), $X^{(1)}_1,\ldots,X_M^{(N)}$ denote $MN$ independent 
indeterminates over {\bf C}). These maps will usually be treated
implicitly as identifications.
With the identification (1.2.2a) we then define, for $a_1,\ldots,a_N$ 
any natural numbers, the 
restriction of $D_{i,j}$ to $S^{a_1}\otimes\ldots
\otimes S^{a_N}$ to be given by
$$D_{i,j}=\sum_{k=1}^M X_k^{(i)}{ \partial\over\partial X_k^{(j)}} :
S^{a_1}\otimes \ldots \otimes S^{a_N} \to S^{a'_1}\otimes \ldots \otimes
S^{a'_N} \eqno (1.2.3) $$
where we have set $$a'_l=\cases{a_i+1, &if $l=i$ \cr
a_j-1, &if $l=j$ \cr a_l &if $l\ne i$ or $j$ . \cr}$$
\smallskip
 (We set $D_{ij}|S^{a_1}\otimes\cdots\otimes S^{a_N} =0$ if $a_j=0$.)
Note that while the operations $X_k^{(i)}$ and  
$\partial\over\partial X_k^{(j)}$ on $S^{\otimes N}$ of course depend
on $\cal B$, the combination (1.2.3) 
is readily verified to be independent of the choice of
basis $\cal B$. Here is a basis-free equivalent definition, which does 
not require that V be finite-dimensional:

Let V be an arbitrary module over an arbitrary commutative ring R.
 
Suppose first that $i<j$; then the natural transformation $D_{i,j}$ may be
equivalently defined (combinatorially rather than in 
the explicit form of a differential operator) as the operator 
$$D_{i,j}:S^{a_1}{\otimes}_R\cdots{\otimes}_R S^{a_N}\to 
S^{a'_1}{\otimes}_R\cdots{\otimes}_R S^{a'_N} $$
which maps $$\sigma^{(1)}\otimes\cdots\otimes\sigma^{(N)}=(v_1^{(1)}\cdot
\ldots\cdot v_{a_1}^{(1)})\otimes(v_1^{(2)}\cdot\ldots\cdot v_{a_2}^{(2)})
\otimes(v_1^{(N)}\cdot\ldots\cdot v_{a_N}^{(N)})\eqno (1.2.4)$$
(with all $v_q^{(p)}\in V$) into the element$$\sum_{\lambda=1}^{a_j}
\sigma^{(1)}\otimes\cdots\otimes(v_{\lambda}^{(j)}\cdot\sigma^{(i)})
\otimes\cdots\otimes (v_1^{(j)}\cdot\ldots\widehat{v_{\lambda}^{(j)}}
\cdot\ldots\cdot v_{a_j}^{(j)})
\otimes\cdots\otimes\sigma^{(N)}$$
with a precisely similar definition (except that the order of the two 
main parentheses is 
reversed) if $i>j$---while if $i=j$, $D_{i,i}$ acts 
on $S^{a_1}\otimes\cdots\otimes S^{a_N}$ as multiplication by $a_i$ .

In all three cases, $D_{i,j}$ operates on an element $$\omega=(v_1^{(1)}\cdot
\ldots\cdot v_{a_1}^{(1)})\otimes\cdots\otimes(v_1^{(N)}\cdot\ldots\cdot v_{a_N}^{(N)})$$ like this:

In all possible ways remove a factor $v_{\lambda}^{(j)}$, where $1\le \lambda \le a_j$,  from the j-th tensorand $\sigma^{(j)}$ in $\omega$, and 
re-insert this 
$v_{\lambda}^{(j)}$ into the $i$-th tensorand; then add all the $a_j$ results thus obtained.
\medskip
{\bf EXAMPLE:}
\smallskip
$D_{13}(x_1x_2x_3\otimes y_1\otimes z_1z_2)=x_1x_2x_3z_1
\otimes y_1 \otimes z_2 + x_1x_2x_3z_2 \otimes y_1 \otimes z_1 $

and \smallskip
$D_{11}(x_1x_2x_3\otimes y_1 \otimes z_1z_2) = 3 x_1x_2x_3\otimes y_1 \otimes z_1z_2$
\par\bigskip
(Note: Some authors prefer to say the same thing in still a third 
more high-falutin' way,
 by expressing these elementary polarizations (in the obvious way)
in terms of the component
$$S^aV\to S^{a-1}V\otimes V,\,x_1\otimes\cdots\otimes x_a\mapsto
\sum_{i=1}^a (x_1\otimes\cdots\hat x_i\cdots x_a)\otimes x_i $$
of the comultiplication map of the Hopf algebra $SV$.)

We assume,for the rest of this section, that the ground-ring is ${\Bbb C}$.
Using any of the preceding  equivalent definitions for $D_{i,j}$, it 
is readily verified that 
these elementary polarization operators on $S^{\otimes N} V$ obey the 
commutation relations $$[D_{ij},D_{kl}]=\delta_{jk}D_{il}-\delta_{il}D_{kj}$$
which imply the existence of an action $P$ of ${\frak A}_N $ on $S^{\otimes N}$
by natural transformations, uniquely specified by $$P(E_{ij})=D_{ij}\,.$$
In fact, we thus obtain, it is well-known,
a natural Lie-algebra monomorphism (here to be referred to as the
{\bf Capelli injection})                         
$$C:{\frak A}_N \to \hbox{\rm Nat Tsf}(S^{\otimes N},S^{\otimes N}),
E_{i,j} \mapsto D_{i,j}\,.\eqno(1.2.5)$$  
\medskip
We shall, for the remainder of this paper, treat the Capelli injection
as an identification. {\bf In particular, this paper will always identify
 $D_{i,j}=P(E_{i,j})$ 
(the natural transformation on $S^{\otimes N}$) with $E_{i,j}\in {\frak A}_N$.}
(Usually, the notation $E_{i,j}$ will be used for both.)

In the next sub-section, the Weyl multi-polarization will be defined by 
specifying its image (via $C$) as a natural endomorphism of $S^{\otimes N}$.
\medskip

{\bf CAUTION:} For fixed V, the action of ${\frak A}_N$ on 
$S^{\otimes N} V$ by $\Bbb C$-linear transformations, need {\it not}
be faithful; it is for this reason that we shall instead use the faithful 
action by natural endomorphisms of the functor $S^{\otimes N}$ (
i.e., work in terms of `generic' V.)   

\bigskip
\centerline{\bf \S1.3 Multi-polarization Operators on $S^{\otimes N}$}
\medskip
We now turn to certain somewhat intricate 
operations on $S^{\otimes N}$, for whose
construction
the earliest source known to the present author is H.Weyl (cf.[Weyl,p.39]).
\bigskip
Let $\Pi_{\pm}^N$) denote the free Abelian group (written additively)
 on the set of $N^2$ matrices $E_{i,j}$ explained in \S1.1. Thus
the elements in $\Pi_{\pm}^N$ are of the form
$$\sigma=\sum_{i,j\in\underline N}\sigma_{i,j}E_{i,j} 
\,\,\hbox{\quad (\rm all } \sigma_{i,j}\in{\Bbb Z}\,); \eqno(1.3.1)$$
each such $\sigma$ may also be regarded as an $N\times N$ 
matrix $||\sigma_{i,j}||$ with integer entries.

Those $\sigma$ for which all $\sigma_{i,j}\ge 0$ will be called
 {\bf N-shifts}; the set of N-shifts will be denoted by $\Pi^N$. (Thus,
 $\Pi^N$ is the free Abelian semi-group on the set of $E_{i,j}$'s.)
An element $\sigma$ in $\Pi_{\pm}^N$ will be called {\bf effective} if
all $\sigma_{i,j}$ are non-negative (i.e., if $\sigma$ is an $N$-shift)
 and {\bf non-effective} otherwise.

In [Weyl, loc.cit.], Weyl associates (in a slightly different notation) to 
every $N$-shift $\sigma \in \Pi^N$, a transformation 
$$P(\sigma):S^{\otimes N}V \to S^{\otimes N}V$$
which is now to be defined, and which we shall call the {\bf Weyl
polarization operator associated to $\sigma$; }we shall also sometimes refer to
these as {\bf multi-polarization operators.} They include, as a special 
case, the elementary polarization operators discussed in the preceding 
section.

It will be convenient, for the $N$-shift given by (1.3.1), to set
$$ \sigma !=\prod _{i,j \in \underline N} (\sigma_{i,j}!) \eqno(1.3.2)$$
(and to define $\sigma !$ to be 0 if $\sigma$ is non-effective.)      

Again, (as in $\S1.2$) we shall give two definitions (equivalent 
where both are defined).
 The first, that given in [Weyl, loc.cit.]; cf. also (Fulton
and Harris[FH, pp.514--515]) requires that V be a
finite-dimensional complex vector-space,and involves a choice of basis
for $V$.The second is defined for $V$
an arbitrary module over an arbitrary commutative ring.
\medskip
 {\bf DEFINITION 1.3.1}

Consider first the case that $V$ is a finite-dimensional
complex vectorspace, with ${\Bbb C}$-basis
$${\cal B}=(x_1,\ldots,x_M)$$
(so $M=$dim $V$). As before, we use ${\cal B}$ to define the
endomorphisms (1.2.2) and (1.2.2a). Now write (in the commutative semi-group
$\Pi^N$)
$$\sigma =E_{i_1,j_1}+E_{i_2,j_2}+\cdots +E_{i_L,j_L}\hbox{.}$$
Then we define the {\bf non-normalized Weyl polarization}
$$P_0(\sigma):S^{\otimes N}V \to S^{\otimes N}V$$
to be the endomorphism of $S^{\otimes N}V$ given by:
$$P_0(\sigma )=\sum_{k_1=1}^M\sum_{k_2=1}^M\cdots\
\sum_{k_L=1}^M X_{k_1}^{(i_1)}X_{k_2}^{(i_2)}\cdots X_{k_L}^{(i_L)}
{\partial \over \partial X_{k_1}^{(j_1)}}{\partial \over \partial
X_{k_2}^{(j_2)}}\ldots{\partial \over \partial X_{k_L}^{(j_L)}}
 \eqno (1.3.3)$$
There is (as will become apparent) some  advantage to considering instead 
the {\bf normalized Weyl polarizations}
$$P(\sigma)\buildrel \rm def\over = ({1 \over \sigma !})P_0(\sigma)\eqno(1.3.4)
$$
(In particular, it is with this normalization that Th.2.3.1 below becomes
valid.)

When mention is made below simply of Weyl
 polarizations, these
normalized polarizations $P(\sigma)$ are always to be understood.
Although $P_0(\sigma)$ and $P(\sigma)$ are basis-independent and natural
in $V$, these facts are perhaps not obvious at this stage (but
will be immediate corollaries of Th. 2.3.1 below).
\smallskip
{\leftline\bf REMARKS ON NOTATION:} Weyl in [Weyl, loc.cit.] uses the notation
$$\Delta_{i_1,j_1}\Delta_{i_2,j_2}\cdots\Delta_{i_N,j_N}$$
for the expression (1.3.3), which he calls a {\bf quasi-composition},
 remarking that (unlike the ordinary composition $D_{i_1,j_1}
\cdots D_{i_N,j_N}$) it is unchanged by rearrangement of the factors.
It seems perhaps more lucid to use instead a notation such as $P_0(\sigma)$
for (1.3.3), particularly since it is {\it not} always true that 
$$P_0(E_{i_1,j_1}+E_{i_2,j_2}) (=\Delta_{i_1,j_1}\Delta_{i_2,j_2}
\hbox{\thinspace in Weyl's notation)}$$
is the composite of $D_{i_1,j_1}=P_0(E_{i_1,j_1})$ 
and  $D_{i_2,j_2}=P_0(E_{i_2,j_2})$.

Weyl's original notation (which was elegantly suited to his proof 
in [Weyl,Chapter II, \S4] of the Capelli identity)
seems in the spirit of the 'umbral' notation
of nineteenth century invariant theory, in which 'symbolical products' occur,
i.e. objects which look like products but are not, and require
re-interpretation.

Also, the additive notation we have adopted in $\Pi^N$, has the advantage of
avoiding confusion between e.g. the $N$-shift
$$E_{1,2}+E_{2,3}=E_{2,3}+E_{1,2}\in\Pi^N$$
and the elements
$$E_{1,2}E_{2,3}\ne E_{2,3}E_{1,2}$$
in the enveloping algebra ${\frak A}_N$. 
\bigskip
\bigskip
By the {\it weight vector} $wt( \sigma )$ of an $N$-shift $\sigma$, will
be meant the $N$-tuple
$$wt(\sigma)=(wt_1(\sigma,\ldots,wt_N(\sigma) \eqno(1.3.5)$$
where
$$\eqalignno{
wt_i(\sigma)&=\sum_{j=1}^N\sigma_{i,j} -\sum_{j=1}^N\sigma_{j,i}
&(1.3.5a) \cr
&=(i^{th} \hbox{ row-sum of } \sigma \hbox{, minus the } i^{th}
\hbox{ column-sum }) \cr
}$$
For $a_1,a_2,\cdots,a_N $ any natural numbers, it is readily verified that 
$P(\sigma)$ maps 
$$S^{a_1}\otimes\cdots\otimes S^{a_N} \to S^{a_1+wt_1(\sigma)}\otimes
\cdots\otimes S^{a_N+wt_N(\sigma)}\eqno(1.3.6)$$
(This is 0 if any $a_i+wt_i(\sigma)=b_i$ is negative, in accordance 
with the usual convention that $S^bV=0$ if $b<0$.)
\smallskip
Note that $\sum_{i=1}^N wt_i(\sigma)=0$.     
\bigskip
Now let us drop the assumption that V is finite-dimensional. Our next goal
is to give a second definition of $P(\sigma)$ using combinatorial 
concepts rather than partial derivatives (but agreeing, as will be proved
below, with Def. 1.3.1 where the latter has been defined.)
\bigskip

\proclaim Definition 1.3.2. If V is a module over the commutative ring R, and 
if $$\alpha =(a_1,\ldots,
a_N)$$ is any $N$-tuple of non-negative integers, then  
 in order to define the restriction
$$P(\sigma,\alpha):=P(\sigma)|(S^{a_1}V\otimes
\cdots\otimes S^{a_N}V)$$
its action is to be given on the generating elements
$$\omega=(x_1^{(1)}\cdot x_2^{(1)}\cdot\ldots\cdot x_{a_1}^{(1)})\otimes\cdots
\otimes(x_1^{(N)}\cdot x_2^{(N)}\cdot\ldots\cdot x_{a_N}^{(N)})$$
\leftline{(with all $x_p^{(q)}$ elements of $V$), by the following rule:} 
Namely: $P(\sigma,\alpha)$ acts on $\omega$ by moving 
(in all possible ways) an {\bf unordered} collection of $\sigma_{i,j}$
letters $x$ from the $j$-th tensor factor of $\omega$ to 
the $i$-th, for all $i$ and $j$ between $1$ and $N$ ---
NO LETTER BEING MOVED TWICE. The results
are then summed, to obtain $P(\sigma)\omega$.

\bigskip
(For the time being, we shall refer to the endomorphisms given by
Definition 1.3.1 as {\bf differential polarizations}, and those given by 
Definition 1.3.2 as {\bf combinatorial polarizations}.This distinction is only
a temporary one, since
in Section 2.3 these two definitions will be proved equivalent
where both are defined. Also, the preceding somewhat breezy version of
Definition 1.3.2 will be restated more formally in Section 1.4.)

Note that (1.3.6) is clearly still the type of $P(\sigma)$, if 
$P(\sigma)$ is interpreted in
the sense of Def.1.3.2.

Perhaps it will be helpful at this point to give some illustrative 
examples for Definition 1.3.2:
\medskip

{\bf EXAMPLE 1.3.1}
\noindent\smallskip
If $$\sigma_1=E_{1,2}+2E_{13}+3E_{32} \in \Pi^3 \hbox{,}$$
i.e. if $$\sigma_1=\pmatrix{
0&1&2 \cr
0&0&0 \cr
0&3&0 \cr
} \hbox{,}$$then $\sigma_1$ has for weight-vector
$$wt(\sigma_1)=(3-0,0-4,3-2)=(3,-4,1)$$
and the combinatorial polarization
$$P(\sigma_1,(a_1,a_2,a_3)):S^{a_1}V\otimes S^{a_2}V\otimes S^{a_N}V \to
S^{a_1+3} \otimes S^{a_2-4}\otimes S^{a_3+1}$$
is, (in accordance with Definition 1.3.2), the map which takes 
$$\omega=(x_1\cdot\ldots\cdot x_{a_1})\otimes(y_1\cdot\ldots\cdot y_{a_2})
\otimes(z_1\cdot\ldots\cdot z_{a_3})$$
(all $x$'s, $y$'s and $z$'s in V) into the sum of all 
$a_2\cdot{a_3 \choose 2}\cdot{a_2-1 \choose 3}$                              
terms obtained by moving:
\smallskip
\noindent one letter y from the second tensorand into the 
first (corresponding 
to the entry
$(\sigma_1)_{1,2}=1$), two letters z from the third tensorand into the first
(corresponding to $(\sigma_1)_{1,3}=2$), and three letters y from the second  
tensorand to the third(corresponding to $(\sigma_1)_{3,2 }=3$), no letter being moved twice.

In other words, $P(\sigma_1)\omega$ is to equal the sum
$$\sum  
\left[(y_{j_1}z_{k_1}z_{k_2}x_1\cdots x_{a_1})\otimes({y_1\cdots y_{a_2}
\over y_{j_1}y_{J_1}y_{J_2}y_{J_3}})\otimes(y_{J_1}y_{J_2}y_{J_3}z_1\cdots
\widehat{z_{k_1}}\cdots\widehat{z_{k_2}}\cdots z_{a_3})\right]\,,
$$ extended over the indexing set indicated by
$$\sum_{1\le{j_1}\le {a_2}}
\sum_{\scriptstyle 1\le{J_1}<{J_2}<{J_3}\le{a_2}
\atop \scriptstyle {j_1}\notin\{J_1,J_2,J_3\}}\sum_{1\le{k_1}<{k_2}\le{a_3}}$$
\smallskip
{\bf REMARKS:}
\smallskip
A)We employ in this formula, (for ease of notation), two quite 
equivalent methods for denoting the deletion of factors from a product:
in the second tensorand deletion of $y_{j_1}y_{J_1}y_{J_2}y_{J_3}$ is
indicated as a division, while in the third tensorand, deletion of
$z_{k_1}z_{k_2}$ is indicated by the usual ``$\wedge$'' notation. We could just
as well have written these the other way around.

B)If $a_2<4$ then this sum is empty, and in accordance with the usual
conventions for an empty sum, $P(\sigma_1)\omega=0$.
 \medskip
{\bf EXAMPLE 1.3.2}
\smallskip \noindent If $$\sigma_2=\pmatrix{2&1 \cr 1&0 \cr}\hbox
{ and }\alpha=(a_1,a_2)$$
then $wt(\sigma_2)=(0,0)$,and $$P(\sigma_2,\alpha):
S^{a_1}V\otimes S^{a_2}V\to S^{a_1}V\otimes S^{a_2}$$ is 
the map which sends
$$\omega=(x_1\cdot\ldots\cdot x_{a_1})\otimes(y_1\cdot\ldots\cdot y_{a_2})$$ into the sum 
$$\sum_{1\le I_1<I_2\le a_1}\sum_{\scriptstyle 1\le i\le a_1\atop
\scriptstyle i\notin\{I_1,I_2\}}\sum_{1\le j\le a_2}
(y_j x_1 \cdots \widehat {x_i}\cdots x_{a_1})\otimes
(x_iy_1\cdots\widehat{y_j}\cdots y_{a_2})
$$ 

Note the effect of the diagonal entry $(\sigma_2)_{11}=2$
in $\sigma_2$, is to select (in all possible ways), an unordered pair
of letters $\{x_{I_1},x_{I_2}\}$ in the first tensorand (both $I_1$
and $I_2$ being distinct from $i_1$), whereupon we proceed to leave these 
two letters untouched. (We may as well suppose $I_1<I_2$.)
Thus $P(\sigma_2)$ differs simply by the scalar factor 
${a_1 -1}\choose 2$ from $P(\sigma_3)$ where
$$\sigma_3=\pmatrix{0&1 \cr 1&0\cr}$$---i.e.,
\smallskip
$$P(\sigma_2,\alpha)={{a_1 -1}\choose 2}P(\sigma_3,\alpha)$$
where
$$P(\sigma_3,\alpha)\omega=\sum_{1\le i\le a_1}\sum_{1\le j\le a_2}
(y_j x_1 \cdots \widehat {x_i}\cdots x_{a_1})\otimes
(x_iy_1\cdots \widehat{y_j}\cdots y_{a_2})\hbox{ . }
$$

\def\s#1#2{\sigma_{#1,#2}}
\def\x#1#2{x^{(#1)}(#2)}
\def\dfeq{\buildrel \rm def \over =}
\def\E#1#2{C_{#1,#2}^{\iota}}
\long\def\ignore#1{}
\bigskip
\centerline{\bf \S 1.4 \ $\sigma$-{\bf selections}}
\medskip

We next need to restate Definition 1.3.2 somewhat more 
formally. For this purpose, and also for later computations, the following 
notations will be convenient:

Let $R$ be a fixed commutative ground-ring, and let $V$ be an arbitrary
$R$-module.
Recall  from $\S 1.1$ the notation
$$\underline N=\{ 1,...,N \}$$
For any map $$x:\underline a \to V$$ we define                    
$$x[[\,\underline a\,]] \buildrel \rm def \over = x(1)\cdot x(2)\cdot\ldots\cdot 
x(a) \in S^a V$$
Note that such elements span $S^a V$ over $R$; and that similarly
(for all non-negative integers $a_1,\cdots,a_N$)
$$S^{ a_1}V\otimes\cdots\otimes S^{a_N}V$$ 
is spanned over $R$ by the set of all elements of the form
$$\eqalignno{
\omega&=x^{(1)}[[\,\underline {a_1}\,]]\otimes\cdots\otimes x^{(N)}
[[\,\underline{a_N}\,]]& (1.4.1)\cr
&=(\x1 1\cdot\ldots\cdot\x1 {a_1})\otimes\cdots\otimes
(\x{N} 1\cdot\ldots\cdot\x{N} {a_N}) \cr}  $$
(where,  for $1\le i\le N$, $x^{(i)}$ is an arbitrary map 
$\underline a_i \to V$).

More generally, for any sub-set $E$ of $\underline a$, say 
$$E=\{e_1,\cdots,e_L\}\subseteq \underline a\,,\hbox{
\quad with $e_1,\cdots,e_L$ distinct,}$$
and any map  $$x:\underline a \to V\,,$$ we define
$$x[[E]] \buildrel \rm def \over = x(e_1)\cdot x(e_2)\cdot\ldots\cdot 
x(e_L) \in S^L V$$

\proclaim Definition 1.4.1. Let
$$\sigma \in \Pi^N,\alpha=(a_1,a_2,\cdots,a_N) \in {\Bbb N}^N\,;$$  
then by a $\sigma$-{\bf selection from} $\alpha$ will be meant
a collection $$\{C_{i,j}:i \hbox{ and }j \in \underline N \}$$
such that, (for all $i$ between 1 and $N$),
$C_{1,i},C_{2,i},\ldots,C_{N,i}$ are pairwise disjoint subsets of
$\underline a_i$ , and such that  (for all $i$ and $j$ between 1 and $N$),
$\sigma_{i,j}$ is the cardinality of $C_{i,j}$. 
We denote by $${\alpha \choose \sigma}$$
the set of all such.

{\bf NOTE:} Given such a $\sigma$-selection
$$\iota =\{\E{i}j :i,j\in \underline N\}\in {\alpha \choose \sigma}
 \eqno(1.4.2)$$
it will be convenient to set, for $1\le i\le N$,
$$\E0i \dfeq \{1,\ldots,a_i\}-( \bigcup_{j=1}^N \E{j}{i}) \eqno(1.4.3)$$
and
$$\sigma_{0,i}\dfeq a_i-\sum_{j=1}^N \sigma_{j,i} \eqno(1.4.4)$$
Note  that thus each $\underline a_i$ is the disjoint union of
$$\E0i\,,\E1i\,,\cdots,\E{N}i\,,$$
and that 
$$\#(\E0i)=\sigma_{0,i}$$
The cardinality of $${\alpha \choose \sigma}$$ is then given by
$$\#{\alpha \choose \sigma}=\prod_{i=1}^N {a_i! \over {\sigma_{1,i}!\cdot\ldots\cdot
\sigma_{N,i}!\cdot(a_i-\sigma_{1,i}-\cdots-\sigma_{N,i})!}}$$
\medskip
Let us now express $P(\sigma)\omega$ (where $\omega$ is given by 
eqn.(1.4.1)), in terms of the notation just explained. Recall that Def.1.3.2 
spoke of ``moving (in all possible ways) $\cdots\sigma_{i,j}$
letters $x$ from the $j$-th tensor factor of $\omega$ to 
the $i$-th, for all $i$ and $j$ between $1$ and $N$ ---
NO LETTER BEING MOVED TWICE.  The results
are then summed, to obtain $P(\sigma)\omega$.''

Clearly, such a `possible way' of moving the letters of
$\omega$, is furnished precisely by a $\sigma$-selection
$$\iota=\E{i}{j} \in {\alpha \choose \sigma}$$  
which instructs us to move the $\sigma_{i,j}$ letters
$$\x{j}{e} \hbox{ where } e \in \E{i}{j}$$
from the $j$-th tensor factor of $\omega$ to 
the $i$-th, thus obtaining a result which we shall denote by
$\omega^{\iota}$, and whose precise value is given by:
$$ \omega^{\iota} \dfeq (\omega^{\iota})_1\otimes\cdots
\otimes(\omega^{\iota})_N \eqno(1.4.5)$$
where (for all $i$ between 1 and $N$), we set 
$$(\omega^{\iota})_i \dfeq (\prod_{j=1}^N x^{(j)}[[\E{i}{j}]])
\cdot x^{(i)}[[\E{0}{i}]] \in 
S^{a_i+wt_i(\sigma)}V \eqno(1.4.6)$$
(cf. eqn.(1.3.5a) in $\S1.3$).
With the notation thus defined, the following equation may then be regarded 
as the deluxe version of Definition (1.3.2):
$$P(\sigma)\omega \dfeq \sum_{\iota \in {\alpha \choose \sigma}}
\omega^{\iota}\,. \eqno(1.4.7)$$
\bigskip
\bigskip
{\bf NOTE:} Our formulas sometimes contain expressions of the form
$P(\sigma')$ where $\sigma'$ is an $N\times N$ matrix over ${\Bbb Z}$, not all
of whose entries are $\ge0$. (In the terminology introduced in the beginning of
$\S1.3$, $\sigma'$ is a non-effective element of $\Pi_{\pm}^N$)

To avoid any possible ambiguity, let us agree always to set, in such a case, 
$$P(\sigma')=0 \hbox{ if any }\sigma'_{i,j}<0\eqno(1.4.8)$$

\def\s#1#2{\sigma_{#1,#2}}
\def\x#1#2{x^{(#1)}(#2)}
\def\dfeq{\buildrel \rm def \over =}
\def\E#1#2{C_{#1,#2}^{\iota}}
\def\sred{\sigma^{\rm red}}
\bigskip
\centerline{\bf Section 2  The Product of a 
Multi-polarization by an Elementary Polarization}
\medskip
\centerline{\bf \S 2.1 The Recurrence for Combinatorial Polarizations}
\medskip
 {\it NOTE:Throughout this subsection,
 all multi-polarizations $P(\sigma)$ 
considered are to be understood in the sense of Definition 1.3.2 of the 
preceding \S1.3.}
 (We shall consider in the following subsection 2.2,
multi-polarizations in the sense of Definition 1.3.1, i.e. in
the sense of Weyl's original definition.)

\proclaim Theorem 2.1.1.  Let $\sigma \in \Pi^N$ and let $i$,$j$ be 
distinct integers between $1$ and $N$.Then 
$$E_{i,j}P(\sigma) =
(\sigma_{i,j}+1)P(E_{i,j}+\sigma)+
\sum_{k=1}^N (\sigma_{i,k}+1)P(\sigma+
E_{i,k}-E_{j,k})  \eqno(2.1.1) $$

{\bf PROOF:} By an obvious symmetry of the situation under study (namely, via
the action of the symmetric group on N letters, upon the tensor product of 
$N$ copies of $SE$), it is evident that it suffices to prove (2.1.1) in
the special case $i=2,j=1$.

Thus, it suffices to verify that
$$E_{2,1}P(\sigma)\omega =
(\sigma_{2,1}+1)P(E_{2,1}+\sigma)\omega+
\sum_{k=1}^N (\sigma_{2,k}+1)P(\sigma+
E_{2,k}-E_{1,k})\omega  \eqno(2.1.2) $$
for every generating element $\omega$ of the form (1.4.1). 
Let us  then consider (utilizing the notation explained in \S1.4)
$$E_{2,1}P(\sigma)\omega=\sum_{\iota \in {\alpha \choose \sigma}}
E_{2,1}\omega^\iota=\sum_{\iota \in {\alpha \choose \sigma}}
E_{2,1}[(\omega^{\iota})_1\otimes
(\omega^{\iota})_2\otimes\cdots\otimes(\omega^{\iota})_N] $$
$E_{2,1}\omega^{\iota}$  is  
 the sum of all terms
obtained by removing one letter u from the first tensorand
$$\eqalignno{
(\omega^\iota)_1&=x^{(1)}[[\,\underline {a_1}-\E21-\cdots-
\E{N}1]]\cdot x^{(2)}[[\E12]]\cdot\ldots\cdot x^{(N)}[[\E1N]] &(2.1.3) \cr
&=x^{(1)}[[\,\underline {a_1}-\E11-\cdots-
\E{N}1]]\cdot \prod_{k=1}^N x^{(k)}[[\E1k]]   \cr }
$$
of $\omega^{\iota}$, and inserting it instead into the second
tensorand
$$(\omega^\iota)_2=x^{(1)}[[\E21]]\cdot x^{(2)}[[\,\underline{a_2}
-\E12-\E32-\cdots-\E{N}2]]\cdot\prod_{k=3}^N x^{(k)}[[\E2k]]\eqno(2.1.3a)$$
---let us denote by  $T(\omega^{\iota},u)$ the term thus obtained.
That is to say, for each letter u in (2.1.3) we set
$$T(\omega^{\iota},u)=({(\omega^{\iota})_1 \over u})\otimes 
(u\cdot(\omega^{\iota})_2)\otimes(\omega^{\iota})_3
\cdots\otimes(\omega^{\iota})_N \eqno(2.1.4)$$
and we then have
$$E_{2,1}P(\sigma)\omega=\sum_{\iota \in {\alpha \choose \sigma}}
\sum_u T(\omega^{\iota},u) \eqno(2.1.5)$$
(the inner sum being extended over all letters u in (2.1.3).)

The next step is to decompose (2.1.5) into $N+1$ sub-sums, (according
to where u comes from), as follows:

If $1\le k\le N$, let us define ${\cal S}(k)$ 
to be the sum of those terms in (2.1.5) for which $u=x^{(k)}(e)$ with
$e\in \E1k$, i.e. we set
$${\cal S}(k)\dfeq \sum_{\iota \in {\alpha \choose \sigma}}
\sum_{e^{(k)} \in \E1{k}} T(\omega^{\iota},x^{(k)}(e^{(k)})) \eqno(2.1.6a)$$
(Here we adopt, as always, the convention which sets empty sums equal 
to 0--- thus ${\cal S}(k)$ is to be 0 if $\alpha\choose\sigma$ is empty,
or if $\E1k=\emptyset$.) 
Similarly, we set 
$${\cal S}(0)\dfeq \sum_{\iota \in {\alpha \choose \sigma}}
\sum_{e^{(0)} \in \E01} T(\omega^{\iota},x^{(1)}(e^{(0)}))
\eqno(2.1.6b)$$
Thus we have the desired decomposition:
$$E_{2,1}P(\sigma)\omega=\sum_{k=o}^N {\cal S}(k) \eqno(2.1.7)$$   

To complete the proof of the theorem, it thus suffices to prove the two 
following equations:
$${\cal S}(0)=(\sigma_{2,1}+1)P(\sigma+E_{2,1})\omega \eqno (2.1.8)$$
and (for $1\le k\le N$)
$${\cal S}(k)=(\sigma_{2,k}+1)P(\sigma+E_{2,k}-E_{1,k})\omega \eqno(2.1.9)$$
\medskip

First, suppose  $1\le k\le N$, and consider ${\cal S}(k)$. In 
proving (2.1.9), we must consider separately two cases, according to
whether or not $\sigma+E_{2,k}-E_{1,k}$ has all entries non-negative. 
\smallskip
\leftline{\bf CASE ONE:\qquad $\underline{\sigma+E_{2,k}-E_{1,k}\in \Pi^N}$ }
\smallskip  
In order to establish (2.1.9),
 we next construct a 
$(\sigma_{2,k}+1)$-to-one correspondence, mapping the
$$ [\#{\sigma \choose \alpha}] \cdot \sigma_{1,k} $$
terms         
$$\{T(\omega^{\iota},x^{(k)}(e^{(k)}): \iota \in 
{\alpha \choose \sigma},e^{(k)} \in \E1k\}$$
of  the sum (2.1.6a) which makes up  ${\cal S}(k)$, into the  
$$\#{\sigma+E_{2,k}-E_{1,k} \choose \alpha}$$
terms 
$$ \{\omega^{\kappa}: \kappa 
\in {\sigma+E_{2,k}-E_{1,k} \choose \alpha} \}$$  
whose sum is $P(\sigma+E_{2,k}-E_{1,k})\omega$, in such a way that 
equality holds for each pair of corresponding terms. 

Namely, for every pair 
$$\iota\in{\alpha\choose\sigma},e^{(k)}\in \E1k \eqno(2.1.10)$$
let us define
$$\kappa=\Phi_k(\iota,e^{(k)})=\{C_{i,j}^{\kappa}:i,j \in \underline N\}$$
to be the $(\sigma+E_{2,k}-E_{1,k})$-selection specified by
$$\cases{
C_{2,k}^{\kappa}=\E2k\cup\{e^{(k)}\} & \cr
C_{1,k}^{\kappa}=\E1k-\{e^{(k)}\} & \cr
C_{i,j}^{\kappa}=\E{i}j\,,&  otherwise. \cr }\eqno(2.1.11) $$
(In other words, $\kappa$ is obtained from $\iota$ by moving $e^{(k)}$
from $\E1k$ to $\E2k$.)

Then, for the selection $\kappa$ thus defined, (1.4.6) obviously becomes
$$\left\{ \eqalign{
(\omega^{\kappa})_1& =x^{(1)}[[\,\underline{ a_1}-\{e^{(k)}\}-\E21-\dots
-\E{N}1]]\cdot x^{(2)}[[\E12]]\cdot\ldots\cdot x^{(N)}[[\E1N]] \cr
(\omega^{\kappa})_2& =x^{(1)}[[\{e^{(k)}\}\cup\E21]]\cdot
 x^{(2)}[[\,\underline{a_2}-\E12-\E32-\cdots-\E{N}2]]\cdot
\prod_{l=3}^N x^{(l)}[[\E2l]]  \cr
(\omega^{\kappa})_l& =\E{i}{j}=(\omega^{\iota})_l  \hbox{ if $3\le l\le N$} \cr
}\right. $$
which implies (using (2.1.3), (2.1.3a) and (2.1.4)) that
$$T(\omega^{\iota},x^{(k)}(e^{(k)})=\omega^{\kappa}$$
\smallskip
Given any
$$\kappa \in {\alpha \choose \sigma+E_{2,k}-E_{1,k}}$$
there are precisely $\sigma_{2,k}+1$ pairs (2.1.10) such that  
$$\kappa=\kappa(\iota,e^{(k)})$$
---given by letting $e^{(k)}$ range through the 
$$(\sigma+E_{2,k}-E_{1,k})_{2,k}=\sigma_{2,k}+1$$
elements of $C_{2,k}^{\kappa}$, and then taking, for each such $e^{(k)}$,
the $C^\iota$'s uniquely determined by (2.1.11). Hence we obtain
$$ \eqalign{
{\cal S}(k)& =\sum_{\iota \in {\alpha\choose\sigma}}\sum_{e^{(k)}\in
\E1k}T(\omega^{\iota},x^{(k)}(e^{(k)})) \cr
& =({\sigma}_{2,,k}+1)\sum_{\kappa \in {\alpha \choose \sigma+E_{2,k}-E_{1,k}}}
\omega^{\kappa} \cr
& =(\sigma_{2,k}+1)P(\sigma+E_{2,k}-E_{1,k})\omega \cr
}$$
(as was to be proved.)
\smallskip
\leftline{\bf CASE TWO:\qquad $\underline{\sigma+E_{2,k}-E_{1,k}\not\in 
\Pi^N}$ }
\smallskip
Since $\sigma\in\Pi^N$ by hypothesis, it must be that $\sigma_{1,k}=0$, whence
{\cal S}(k), being an empty sum, is equal to zero.Since also  
$P(\sigma+E_{2,k}-E_{1,k})=0$, (2.1.9) again holds trivially in this case.
\medskip
Thus (2.1.9) holds in all cases.
We have still to prove (2.1.8). The proof is essentialy the same, with
these two variations:\par\noindent
We must replace (2.1.10) by
$$\iota\in{\alpha\choose\sigma},e^{(0)}\in\E10=\underline{a_1}-
\E11-\E12-\ldots-\E1{N} \eqno(2.1.10a)$$
and we must replace (2.1.11) by
$$\hbox{(2.1.11a)}\hskip 15pt \left\{ \eqalign{
C_{2,1}^{\kappa}& =\E21\cup\{e^{(0)}\}  \cr
C_{i,j}^{\kappa}& =\E{i}j \hbox{  otherwise} \cr }\right. $$
noting that for each 
$$\kappa\in{\alpha\choose\sigma+E_{21}}$$
there are precisely $\sigma_{21}+1$ pairs (2.1.10a) such that 
(2.1.11a) holds (given by letting $e^{(0)}$ range through 
the $\sigma_{21}+1$ elements of $C_{2,1}^{\kappa}$.)
\smallskip
\centerline{This completes the proof of the theorem.}
\bigskip
\medskip

\proclaim Definition 2.1.2. An $N$-shift $\sigma$ will be called 
{\bf reduced} if all its diagonal entries $\sigma_{i,i}$ vanish (for 
$1\le i\le N$). The {\bf reduced form} $\sigma^{\rm red}$ of any 
$N$-shift $\sigma$, is the $N$-shift specified by
$$(\sigma^{\rm red})_{i,j}=\cases{\sigma_{i,j},& if $i\ne j$;\cr
                               0& otherwise. \cr}$$

\bigskip
\proclaim Proposition 2.1.3. Let V be a module over a commutative 
ring R. Suppose $a_1,\ldots,a_N$ are non-negative integers,and that
$$\sigma\in\Pi^N,\omega\in\ S^{a_1}V\otimes\cdots\otimes
S^{a_N}V $$
Then 
$$P(\sigma)\omega=\prod_{i=1}^N {a_i-\sum_{j\ne i}\sigma_{j,i}
\choose \sigma_{i,i}}\cdot P(\sigma^{\rm red})\omega \eqno (2.1.12 )$$  

\noindent [Note: For an example of this proposition, cf. Example 1.3.2]

\noindent
{\bf PROOF:} Let us compare  
$$P(\sigma)\omega=\sum_{\iota\in{\alpha\choose\sigma}}
\omega^{\iota},\,\hbox{ with }\,
P(\sigma^{\rm red})\omega=\sum_{\kappa\in{\alpha\choose{\sigma}^{red}}}
\omega^{\kappa}\,. $$
We have the epimorphism
$$\phi:{\alpha\choose\sigma}\to{\alpha\choose\sred}$$
which maps a $\sigma$-selection $\iota$ into the $\sred$-selection
$$\kappa=\phi(\iota)=\{C_{i,j}^{\kappa}\}$$
defined by
$$C_{i,j}^{\kappa}=\cases {C_{i,j}^{\iota},&if $i\neq j$ \cr
\emptyset& if $i=j$. \cr}$$
Clearly $\omega^{\kappa}=\omega^{\iota}$.\hfil\break
Given $$\kappa\in{\alpha\choose\sred}$$
we obtain all $\iota$ with $\phi(\iota)=\kappa$ by
selecting as $C_{i,i}^{\iota}$, for
each $i$ from 1 to $N$, an arbitrary $\sigma_{i,i}$-element subset
of
$$\underline{a_i}-\bigcup_{\scriptstyle {1\le j\le N}\atop \scriptstyle
j\neq i}C_{j,i}^{\kappa}\,.)$$ 

Thus, $\phi$ is an $M$-to-1 map, where
$$M=\prod_{i=1}^N {a_i-\sum_{j\ne i}\sigma_{j,i}
\choose \sigma_{i,i}}\,,$$
and so we have, as asserted,
$$P(\sigma)\omega=\sum_{\iota\in{\alpha\choose\sigma}}\omega^{\iota}
=M\cdot\sum_{\kappa\in{\alpha\choose\sred}}\omega^{\kappa}
=MP(\sred)\omega\,.$$
\bigskip
\proclaim Theorem 2.1.4.  Let $\sigma \in \Pi^N$ and let $1\le i\le N$. 
Then
$$ E_{i,i}P(\sigma) =(\sigma_{i,i}+1)P(E_{i,i}+
\sigma)+(\sum_{k=1}^N \sigma_{i,k})
P(\sigma) \eqno (2.1.13)$$ 

(It is to be noted that the coefficients in (2.1.1) and (2.1.13) do
not depend, as one might have expected, on  $a_1,\cdots,a_N$, but
 only on $\sigma$.)

\medskip
 \leftline{\bf PROOF:}
It suffices, by the same symmetry argument as before, to prove the special case
$$E_{1,1}P(\sigma)\omega=(\sigma_{1,1}+1)P(E_{1,1}+\sigma)\omega
+(\sum_{k=1}^N \sigma_{1,k})P(\sigma)\omega \eqno(2.1.14) $$
(where $\omega$ is given by (1.4.1)).

Let us now set
$$P'\,\buildrel \rm def \over =\,\prod_{i=2}^N {a_i-\sum_{j\ne i}\sigma_{j,i}
\choose \sigma_{i,i}}\cdot P(\sigma^{\rm red})\omega $$  
Note that (by (1.4.4))
$$\underline a_1-\sigma_{2,1}-\cdots-\sigma_{N,1}=\sigma_{1,1}+\sigma_{0,1}$$
Hence, Prop.2.1.3 implies
$$P(\sigma)\omega={a_1-\sigma_{2,1}-\cdots-\sigma_{N,1} \choose
\sigma_{1,1}}P'={\sigma_{1,1}+\sigma_{0,1} \choose \sigma_{1,1}}P'
\eqno(2.1.15)$$
Replacing $\sigma$ by $\sigma+E_{1,1}$ in this last equation yields
$$P(\sigma+E_{1,1})\omega={\sigma_{1,1}+\sigma_{0,1}
 \choose \sigma_{1,1}+1}P' \eqno(2.1.16)$$   

Recall  that $E_{1,1}$ acts on $$ S^{a_1}V \otimes
\cdots\otimes S^{a_N}V \eqno(2.1.17)$$
as multiplication by $a_1$. Since 
$\omega$ lies in (2.1.17), it follows by (1.3.6) (which as noted 
applies both for Def.
1.3.1 and for Def.1.3.2) that
$$P(\sigma)\omega\in S^{a_1+wt_1(\sigma)}V 
\otimes\cdots\otimes S^{a_N+wt_N(\sigma)}V$$  
Here $wt_i(\sigma)$ is defined  by (1.3.5a); we 
note in particular
$$\eqalignno{
a_1+wt_1(\sigma)&=a_1-(\sigma_{2,1}+\cdots+\sigma_{N,1})
+(\sigma_{1,2}+\cdots+\sigma_{1,N}) &(2.1.18)\cr
&= \sigma_{0,1}+\sum^N_{k=1} \sigma_{1,k}   \cr }$$
From the preceding, we readily deduce 
$$E_{1,1}P(\sigma)\omega=[a_1-(\sigma_{21}+\cdots+\sigma_{N,1})
+(\sigma_{1,2}+\cdots+\sigma_{1,N}]{\sigma_{1,1}+\sigma_{0,1}
 \choose \sigma_{1,1}}P' \eqno(2.1.19)$$

                  In the identity
$$(B+1){A \choose B+1}=(A-B){A \choose B}   \eqno(2.1.20) $$
replace $A$ by $\sigma_{11}+\sigma_{01}$, $B$ by $\sigma_{11}$.
We obtain
$$ (\s11+1){\s11+\s01 \choose \s11+1}=\s01{\s11+\s01 \choose \s11} $$
whence (using also eqns.15, 16 and 19:)
$$ \eqalign{
E_{1,1}P(\sigma)\omega&=(\s01+\sum^N_{k=1}\s1k){\s00+\s01 \choose
 \s11}P' \cr
&=[(\s11+1){\s11+\s01 \choose \s11+1}+(\sum_{k=1}^N \s1k)
\s01{\s11+\s01 \choose \s11}]\cdot P' \cr
&=\s11 P(\sigma+E_{1,1})\omega+(\sum^N_{k=1}\s1k))P(\sigma)\omega
\cr }$$
which completes the proof of the theorem.

\def\x#1#2{X_{#1}^{(#2)}}
\def\d#1#2{\partial \x{#1}{#2}}    
\def\D#1#2{{\partial \over \d{#1}{#2}}}

\bigskip
\centerline{\bf \S2.2 The Recurrence for Differential Polarizations}
\medskip
We assume, throughout the present sub-section, that a specific ${\Bbb C}$-basis 
$${\cal B} =\{x_1,\ldots,x_M\}$$
has been chosen for V, with all Weyl polarizations defined via
Def.1.3.1, by means of
${\cal B}$. (Only in the following $\S2.3$ will we prove, using 
the results of the present sub-section, that the differential
Weyl polarizations are basis-independent.)
\smallskip
Actually, for the purposes of the present sub-section, computations go a bit more 
easily with the {\it non-normalized} differential 
Weyl polarizations $P_0(\sigma)$:
\smallskip
\proclaim Theorem 2.2.1. Let $i$ and $j$ be distinct integers between 
1 and $N$; let $\sigma \in \Pi^N$.Then the four following equations hold:
$$E_{i,j}P_0(\sigma)=P_0(E_{i,j}+\sigma)+\sum_{k=1}^N
\sigma_{j,k}P_0(\sigma+E_{i,k}-E_{j,k}) \eqno(2.2A) $$
$$P_0(\sigma)E_{i,j}=P_0(E_{i,j}+\sigma)+\sum_{k=1}^N
\sigma_{k,j}P_0(\sigma+E_{k,j}-E_{k,i}) \eqno(2.2B) $$
$$E_{i,i}P_0(\sigma)=P_0(E_{i,i}+\sigma)+
\left ( \sum_{k=1}^N \sigma_{i,k}\right )P_0(\sigma)
\eqno(2.2C) $$
$$P_0(\sigma)E_{i,i}=P_0(E_{i,i}+\sigma)+
\left ( \sum_{k=1}^N \sigma_{k,j}\right )P_0(\sigma)
\eqno(2.2D) $$

{\bf PROOF:}\hfil\break
Suppose \hfil
$$\sigma=E_{i_1,j_1}+\cdots+E_{i_L,j_L}\,,\eqno(2.2.1)$$ with all
$i$'s and $j$'s in $\underline N$ (so that
$\sigma_{i,j}$ is the number of $\lambda$ between 1 and $L$
for which $(i_\lambda,j_\lambda)=(i,j)$ ).

For every 
$$F \in {\Bbb C}[\x11,\cdots,\x{M}1;\cdots;\x1N,\cdots,\x{M}N ]
=S^{\otimes N}V\,,$$
we have (by eqn.(1.3.3))
$$E_{i,j}P_0(\sigma)F=(\sum_{k=1}^M \x ki\D kj)\sum_{k_1,\cdots,k_L \in
\underline M}\x{k_1}{i_1}\cdots\x{k_L}{i_L}
{\partial^L F \over \d{k_1}{j_1}\cdots\d{k_L}{j_L}}$$
which is in turn equal to the sum $A+B$ of
$$A=\sum_{k,k_1,\cdots,k_L \in
\underline M}\x ki\x{k_1}{i_1}\cdots\x{k_L}{i_L}
{\partial^{L+1}F \over \d kj\d{k_1}{j_1}\cdots\d{k_L}{j_L}}$$
and  
$$B=\sum_{1\le \lambda\le L}\,
\sum_{k_1,\cdots,k_L \in \underline M}
\delta_{i_{\lambda},j}\x{k_\lambda}{i}\cdot\x{k_1}{i_1}
\cdots\widehat{\x{k_\lambda}{i_\lambda}}\cdots\x{k_L}{i_L}  
{\partial^L F \over \d{k_1}{j_1}\cdots\d{k_\lambda}{j_\lambda}
\cdots\d{k_L}{j_L}}$$
Clearly $A=P_0(E_{i,j}+\sigma)F$, while
$$B=\sum_{\scriptstyle i_{\lambda}=j \atop \scriptstyle 1\le \lambda\le L}
P_0(E_{i_1,j_1}+\cdots+
\underbrace{E_{i,j_{\lambda}}}_{\hbox{replacing\ } E_{i_{\lambda},j_{\lambda}}}
+\cdots+E_{i_L,j_L})F$$ 
(in which the selected term $E_{i,j_{\lambda}}$ replaces the term
$E_{i_{\lambda},j_{\lambda}}=E_{j,j_{\lambda}}$ in (2.2.1).)
Thus,
$$E_{i,j}P_0(\sigma)=P_0(E_{i,j}+\sigma)+\sum_{\scriptstyle i_{\lambda}=j \atop \scriptstyle 1\le \lambda\le L}P_0(\sigma+E_{i,j_{\lambda}}-E_{j,j_{\lambda}})
\eqno(2.2.2) $$
Now,given any k between 1 and $N$, there are precisely $\sigma_{j,k}$ values 
of $\lambda$ for which $i_\lambda=j, j_\lambda=k$, and each contributes
the same term $P_0(\sigma+E_{i,k}-E_{j,k})$ to the sum on the right side of
eqn.(2.2.1); from which eqn.(2.2A) is immediate.

The proof of (2.2B) is precisely similar. For the same reasons, it
suffices to examine here only one of (2.2C) and (2.2D); let us choose the 
first. Then with notation as before,
$$E_{i,i}P_0(\sigma)F=(\sum_{k=1}^M \x ki\D ki)\sum_{k_1,\cdots,k_L \in
\underline M}\x{k_1}{i_1}\cdots\x{k_L}{i_L}
{\partial^L F \over \d{k_1}{j_1}\cdots\d{k_L}{j_L}}$$
which is the sum $A'+B'$ of
$$A'=\sum_{k,k_1,\cdots,k_L \in
\underline M}\x ki\x{k_1}{i_1}\cdots\x{k_L}{i_L}
{\partial^{L+1}F \over \d ki\d{k_1}{j_1}\cdots\d{k_L}{j_L}}=P_0(E_{i,i}+
\sigma)F$$
and 
$$
\eqalign{
B'&=\sum_{ 1\le \lambda\le L}\,
\sum_{k_1,\cdots,k_L \in \underline M}\delta_{i_{\lambda},i}
\x{k_\lambda}{i}\x{k_1}{i_1}
\cdots\widehat{\x{k_\lambda}{i_\lambda}}\cdots\x{k_L}{i_L}  
{\partial^L F \over \d{k_1}{j_1}\cdots\d{k_\lambda}{j_\lambda}
\cdots\d{k_L}{j_L}} \cr
&=\sum_{\scriptstyle i_{\lambda}=i \atop \scriptstyle 1\le \lambda\le L}
\sum_{k_1,\cdots,k_L \in \underline M}\x{k_1}{i_1}\cdots
\x{k_L}{i_L}  {\partial^L F \over \d{k_1}{j_1}\cdots\d{k_\lambda}{j_\lambda}
\cdots\d{k_L}{j_L}} \cr
}$$
Thus ,
$$B'=\sum_{\scriptstyle i_{\lambda}=i \atop \scriptstyle 1\le \lambda\le L}
P_0(\sigma)F=\left ( \sum_{k=1}^N \sigma_{i,k}\right )P_0(\sigma)F$$
and so finally,
$$E_{i,i}P_0(\sigma)=P_0(E_{i,i}+\sigma)+(\sum_{k=1}^N 
\sigma_{i,k})P_0(\sigma)$$
This completes the proof of Th.2.2.1.\par

We obtain the analogous formulas for the normalized differential polarizations
$$P(\sigma)={1 \over \sigma !}P_0(\sigma)$$
as a trivial consequence of the preceding; note 
the first and third of the following formulas, correspond 
precisely with the results in \S2.1 for combinatorial polarizations.
(Let us note, however, that the proofs of the corresponding results in
 \S2.1 were rather more intricate than those in the present sub-section.)
 
\proclaim Corollary 2.2.2. { Let $\sigma \in \Pi^N$ and let $i$,$j$ be 
distinct integers between $1$ and $N$.Then the four following equations hold: 
$$\eqalignno{
E_{i,j}P(\sigma) &=
(\sigma_{i,j}+1)P(E_{i,j}+\sigma)+
\sum_{k=1}^N (\sigma_{i,k}+1)P(\sigma+
E_{i,k}-E_{j,k}) & (2.2.2A) \cr
P(\sigma)E_{i,j} &=(\sigma_{i,j}+1)P(E_{i,j}+\sigma)+
\sum_{k=1}^N (\sigma_{k,j}+1)P(\sigma+E_{k,j}-E_{k,i}) & (2.2.2B) \cr
E_{i,i}P(\sigma) &=(\sigma_{i,i}+1)P(E_{i,i}+
\sigma)+(\sum_{k=1}^N \sigma_{i,k})
P(\sigma) & (2.2.2C) \cr
P(\sigma)E_{i,i} &=(\sigma_{i,i}+1)P(E_{i,i}+\sigma)+
(\sum_{k=1}^N \sigma_{k,i})P(\sigma) & (2.2.2D) \cr
}$$
}

{\bf PROOF:}
Dividing both sides of (1.4A) by $\sigma !$, and observing that
$$(\sigma+E_{i,j})!=\prod_{p,q} (\sigma+E_{i,j})_{p,q}=\sigma !(\sigma_{i,j}
+1\,,)$$ and that
$$(\sigma+E_{i,k}-E_{j,k})!=\sigma !\cdot{\sigma_{i,k}+1 \over \sigma_{j,k}}
$$ we immediately obtain (2.2.2A).

The proof of (2.2.2B) is precisely similar. Similarly,
dividing both sides  of (2.2C), (2.2D) by $ \sigma !$ , we obtain 
(2.2.2C),(2.2.2D).
This completes the proof of Cor.2.2.2.

Finally,as an immediate consequence of the preceding corollary, we obtain
the commutators given by:
\proclaim Corollary 2.2.3. Let $i$ and $j$ be distinct integers between 
1 and $N$; let $\sigma \in \Pi^N$.Then
$[E_{p,q},P(\sigma)]=A-B$, where
$$A=\sum_{k=1}^N (\sigma_{p,k}+1)P(\sigma+E_{p,k}-E_{q,k})\,,$$
and
$$B=\sum_{k=1}^N (\sigma_{k,q}+1)P(\sigma+E_{k,q}-E_{k,p})\,,$$
Also, we have
$$[E_{p,p},P(\sigma)]=\left (\sum_{k=1}^N (\sigma_{p,k}-\sigma_{k,q})
\right )P(\sigma) \eqno(2.2.3)$$\par
\medskip 
{\bf REMARK:} Suppose $i\ne j$. If $\sigma_{j,k}=0$, then the term 
$P_0(\sigma +E_{i,k}-E_{j,k})$ occurs with zero coefficient in
(2.2.A), but $P(\sigma +E_{i,k}-E_{j,k})$ occurs with positive coefficient
$\sigma_{i,k}+1$ in (2.2.2A) (which is a scalar multiple of (2.2A).)
To explain this apparent discrepancy, note that here
$(\sigma +E_{i,k}-E_{j,k})_{j,k}=-1$, so by the convention adopted in
\S1.3, $$P(\sigma +E_{i,k}-E_{j,k})=0\,.$$

\def\s#1#2{\sigma_{#1,#2}}
\def\x#1#2{x^{(#1)}(#2)}
\def\dfeq{\buildrel \rm def \over =}

\bigskip
\centerline{\bf \S2.3 Equivalence of the Two Kinds of Multi-polarizations}
\medskip
\proclaim Theorem 2.3.1. Let $\sigma \in \Pi^N$, and let V be a finite-
dimensional complex vector-space, with
${\Bbb C}$-basis
$${\cal B} =\{x_1,\ldots,x_M\}\,.$$
Let 
$$P^{\cal B}(\sigma,V) \hbox { , resp.  } P(\sigma,V) $$
denote the natural transformations 
$$S^{\otimes N} \to S^{\otimes N}$$
constructed respectively in
Definition 1.3.1 (of  differential Weyl polarizations)
and Definition 1.3.2 (of  combinatorial Weyl polarizations).
Then, $$ P^{\cal B}(\sigma,V)=P(\sigma,V) \eqno(2.3.1) $$

\leftline{\bf PROOF:} By the {\bf weight}\  $W(\sigma)$ of an $N$-shift 
$\sigma$, will be meant the sum of all its entries:
$$W(\sigma) \dfeq \sum_{i=1}^N \sum_{j=1}^N \sigma_{i,j} $$  
We now shall prove (2.3.1) by an induction on $W(\sigma)$: \par
\leftline{\bf Assume first: \quad $\underline{W(\sigma)=1}$}
Here there exist $i$ and $j$ in $\underline N$ such that $\sigma=E_{i,j}$,
and so (2.3.1) becomes the assertion that (as already remarked in
\S1.2) the two definitions for the elementary polarizations 
$D_{i,j}=P(E_{i,j})$ coincide whenever they both make sense.
\medskip

\leftline{\bf Assume next : \quad $\underline{ 1<W(\sigma)\hbox {  and
(2.3.1) holds for all $N$-shifts of weight}< W(\sigma)}$}
\smallskip
There exist i,j in $\underline N$ such that $\sigma_{i,j} >0$; we
consider two cases, according as $i\ne j$ or $i=j$.
\smallskip Suppose first that $i$ and $j$ are distinct.We may write
$$\sigma={\sigma}'+E_{i,j}$$
with ${\sigma}'$ effective. Replacing $\sigma$ by ${\sigma}'$ in eqns.
(2.1.1) and (2.2A) we obtain the two equations
$$\openup2pt \displaylines{
\sigma_{i,j}P(\sigma,V)=E_{i,j}P({\sigma}',V)-
\sum_{k=1}^N ({\sigma}'_{i,k}+1)P({\sigma}'+E_{i,k}-E_{j,k},V)\,, \cr
\sigma_{i,j}P^{\cal B}(\sigma,V)=E_{i,j}P^{\cal B}({\sigma}',V)-
\sum_{k=1}^N ({\sigma}'_{i,k}+1)P^{\cal B}({\sigma}'+E_{i,k}-E_{j,k},V) \cr
}$$
All $N$-shifts on the right sides of these equations, have weights one less
than that of $\sigma$, so the induction hypothesis implies these right sides 
are equal. Hence also the left sides are equal. Since $\sigma_{i,j}$ is by
assumption non-zero, (2.3.1) follows in the present case.

Suppose next that $i=j$. Replacing $\sigma$ by the effective $N$-shift 
$${\sigma}'=\sigma-E_{i,i}$$ 
in eqns.(2.1.13) and (2.2C), we obtain the two equations 
$$\openup2pt \displaylines{
\sigma_{i,i}P(\sigma,V)=E_{i,i}P({\sigma}',V)-
\sum_{k=1}^N ({\sigma}'_{i,k})P({\sigma}',V)\,, \cr
\sigma_{i,i}P^{\cal B}(\sigma,V)=E_{i,j}P^{\cal B}({\sigma}',V)-
\sum_{k=1}^N ({\sigma}'_{i,k})P^{\cal B}({\sigma}',V) \cr
}$$
and conclude the argument as before.

This completes the proof of Theorem 2.3.1.
\bigskip
From this point on, it is no longer necessary (over the ground-field
${\Bbb C}$ ) to distinguish between the two types of multi-polarizations.

The following useful corollary is an immediate consequence of the preceding
argument:

\proclaim Corollary 2.3.2. Let $\sigma$ be an $N$-shift of weight $W>1$,
and let $i$,$j$ in $\underline N$ 
be such that $\sigma_{i,j}$ is non-zero. Then there exist:
effective $N$-shifts (all of weights W-1)
$$\sigma_{1}, \sigma_{2},\cdots, \sigma_{L}$$ 
(for some non-negative integer $L$, which may be 0) and positive
integers $m_k \,(1\le k\le L)$ such that
$$\sigma_{i,j}\cdot P(\sigma)=E_{i,j}P(\sigma -E_{i,j})-
\sum_{k=1}^L m_k P(\sigma_k)\,.$$
If, in addition, $\sigma$ has the property
$$ \hbox{for all }i,j \in \underline N\,,\;\sigma_{i,j}\ne 0 \Rightarrow
i>j \,.$$  
then $\sigma_1,\cdots,\sigma_L$ may be chosen, all also to have
this property.\par

This in turn immediately implies:
\proclaim Corollary 2.3.3. For every $N$-shift $\sigma$ there exists
unique $P'(\sigma)$ in ${\frak A}_N$ such that, for every complex
vector-space V, the Weyl polarization
$$P(\sigma,V):S^{\otimes N}V \to S^{\otimes N}V $$
coincides with multiplication by $P'(\sigma)$.\par

\bigskip
\centerline{\bf Chapter 3  APPLICATIONS OF WEYL POLARIZATIONS }
\medskip
\centerline{\S 3.1 \bf DIFFERENTIALS IN THE ZELEVINSKY COMPLEX}
\medskip
As promised in the introduction, we shall now utilize the Weyl 
polarizations, to obtain new and quite explicit formulas, for 
the differentials in
certain  natural complexes which (in Zelevinsky's metaphor, [Zel,p.152])
`materialize' the Jacobi-Trudi identity. This application involves a
train of thought spanning roughly a century and a half, and one which 
has (it would seem) not yet reached its full conclusion, despite relevant 
work by ({\it inter alia}) Akin[Akin1,2], 
Buchsbaum, Doty[Doty and Doty2], Lascoux[Las],
Maliakas[Mal], Nielsen[Nielsen], Verma[Verma1,2], Woodcock[Wood]
and Zelevinsky[Zel]. 
\medskip
Before discussing the {\it differentials} in the Zelevinsky  
complex, we must begin by 
discussing the individual {\it terms}:

Let us rewrite the Jacobi-Trudi identity (1) in the Introduction, as
$$s_{\alpha}=\sum_{\pi \in {\frak S}_N}sgn(\pi)\prod_{i=1}^N
h_{a_i-i+\pi(i)} \eqno(3.1.1)$$
Reading this in the Grothendieck ring  ${GL_{{\Bbb C}}}^{\wedge}$ , 
$s_{\alpha}$ corresponds to the
Weyl functor $V\mapsto V^\alpha$ (in the terminology of Carter and 
Lusztig,[CL])
while the symmetric polynomial
$$\prod_{i=1}^N h_{a_i-i+\pi(i)}=h_{a_1-1+\pi(1)}h_{a_2-2+\pi(2)}
\cdots h_{a_N-N+\pi(N)}$$
occurring in eqn.(3.1.1), is the character of the functor
$$V \mapsto S^{a_1-1+\pi(1)}V \otimes S^{a_2-2+\pi(2)}V\otimes\
\cdots\otimes S^{a_N-N+\pi(N)}V \eqno(3.1.2)$$
(where as usual, $SV=\oplus S^i V$ denotes the symmetric algebra on the
finite-dimensional complex vector-space $V$.)
 It will be convenient to denote
by ${\cal SYM}^\alpha V$, the vector space
$$S^{a_1}V\otimes S^{a_2}V\otimes\cdots\otimes S^{a_N}V \,,$$
so that (3.1.2) may thus be denoted by
${\cal SYM}^{\alpha -\rho +\pi(\rho)}V$---where we set
$$\rho=(1,2,\cdots,N) \hbox{ and }\pi( \rho )=(\pi(1),\pi(2),\cdots,\pi(N))
\,.$$

All this suggests that, (temporarily forgetting about the
differentials), we may hope to obtain a complex
$ZEL({\alpha})=ZEL(\alpha,V)$ in ${GL_{{\Bbb C}}}^{\wedge}$, the
alternating character-sum of whose terms corresponds as desired to the sum
in (3.1.1), if we define
the $i$-th term of this complex to be
$$ZEL^i({\alpha},V)=\bigoplus_{\scriptstyle \pi \in {\frak S}_N
\atop l(\pi)=i}{\cal SYM}^{\alpha -\rho +\pi(\rho)}V \eqno(3.1.3) $$
(where $l(\pi)$ denotes the number of inversions of a permutation $\pi$
in ${\frak S}_N$.) It will be convenient to denote the individual summands
in (3.1.3) by
$$ZEL^{\pi}(\alpha,V) \buildrel \rm def \over {=}
{\cal SYM}^{\alpha -\rho +\pi(\rho)}V=
S^{a_1-1+\pi(1)}V \otimes S^{a_2-2+\pi(2)}V\otimes\
\cdots\otimes S^{a_N-N+\pi(N)}V \eqno(3.1.3a)$$
\bigskip
Thus we are led (all the authors cited above
appear to agree on this)
to seek for differentials $d_k$ which 
will render exact the following complex:
$$ZEL(\alpha,V):\cdots \buildrel d_{k+1} \over \longrightarrow ZEL^k(\alpha,V)
 \buildrel d_k \over \longrightarrow \cdots \buildrel d_2
\over \longrightarrow ZEL^1(\alpha,V) \buildrel d_1 \over
\longrightarrow ZEL^0(\alpha,V) \eqno(3.1.4) $$ 
and for which $d_1$ has cokernel the Weyl module $V^{\alpha}$
Here we also wish to require the differentials $d_k$ to
be natural transformations---something required by all but one of the
authors cited.
\bigskip
There is a further generalization in ([Zel]); namely,
Zelevinsky next drops the requirement that the $N$-tuple 
$$\alpha=(a_1,\cdots,a_N)$$
represents a partition, requiring only that the $a$'s be integers---
the complex (3.1.4) constructed by Zelevinsky remains exact in this
greater generality (while if in addition $\alpha$ is a partition, then
$\hbox{coker}(d_1)=V^{\alpha}$). The remainder of the present
discussion is to be understood in this greater generality.
\bigskip
To specify the differentials 
$$d_l:ZEL^l({\alpha},V)=\bigoplus_{\scriptstyle \pi \in {\frak S}_N
\atop l(\pi)=l}{\cal SYM}^{\alpha -\rho +\pi(\rho)}V \to
ZEL^{l-1}({\alpha},V)=\bigoplus_{\scriptstyle \pi' \in {\frak S}_N
\atop l(\pi')=l-1}{\cal SYM}^{\alpha -\rho +\pi'(\rho)}V$$
is the same, as to specify the collection of {\Bbb C}-linear
transformations
$$\displaylines{
d_{\alpha}^{\pi,\pi'}:{\cal SYM}^{\alpha -\rho +\pi(\rho)}V=
S^{a_1+\pi(1)-1}V\otimes S^{a_2+\pi(2)-2}V\otimes\cdots\otimes
  S^{a_N+\pi(N)-N}V \cr
\rightarrow {\cal SYM}^{\alpha -\rho +\pi'(\rho)}V=
S^{a_1+\pi'(1)-1}V\otimes S^{a_2+\pi'(2)-2}V\otimes\cdots\otimes
  S^{a_N+\pi'(N)-N}V  \cr }$$
for all 
$$\pi,\pi'\in{\frak S}_N\hbox{ with }l(\pi)=l,\,l(\pi')=l-1  \eqno(3.1.5)$$
(We omit $l$ from our notation for these maps, since by (3.1.5) $l$
is determined by $\pi$.)
\smallskip
We divide the construction of the maps $d_{\alpha}^{\pi,\pi'}$ 
into four  smaller parts, as follows:
\item{A.} We must specify, for precisely which pairs $\pi,\pi'$
the map $d_{\alpha}^{\pi,\pi'}$ is to be non-zero.
\item{B.}To each such pair $\pi,\pi'$ meeting condition A, Akin
assigns a signature $\pm$, according to the rule to be described below.
\item{C.} When $\pi,\pi'$ meet condition A, we shall express 
$d_{\alpha}^{\pi,\pi'}$ as a ${\Bbb C}$-linear combination of
Weyl polarizations $P(\sigma)$---thus we must specify precisely
which $P(\sigma)$ appear in this linear combination. This condition 
will be independent of $\alpha$, and
such shift-matrices $\sigma$ will be called 
{\it subordinate} to the pair  $\pi,\pi'$.
\item{D.} Finally, when $P(\sigma)$ has been designated
 as occurring in $d_{\alpha}^{\pi,\pi'}$, we must specify the precise numerical
coefficient with which it occurs. (It will depend on $\alpha$, and
will be an integer if all $a_i$ are integers---indeed, will be a product
of certain binomial coefficients and factorials, as specified below.)
\bigskip
\centerline{\bf A. WHEN IS  $d_{\alpha}^{\pi,\pi'}$ NON-ZERO?}

Let (3.1.5) hold.

In the complex to be constructed here (which will be proved later to
coincide with that constructed
by Zelevinsky ) the partial maps $d_{\alpha}^{\pi,\pi'}$ are
defined to be 0 unless $\pi'$ precedes $\pi$ in the Bruhat order 
(This {\it Ansatz} was 
suggested to the author by Verma in a conversation.) 

It is well known (cf., for example, [Mathas,p.2,Prop.1.3]) 
that this holds if and
only if there exist i,j such that
 $$1\le i<j\le N,\pi(i)>\pi(j),\,\pi'=\pi(i,j) \eqno(3.1.6)$$

We shall call $(\pi,\pi')$ an {\bf arrow-pair} if these conditions
(3.1.5), (3.1.6)  are satisfied. (The relation of such arrow-pairs
to the theory of maps between Verma modules of $A_N$ will be discussed below
in \S3.4)
\medskip
\centerline{\bf B. AKIN'S NORMALIZATION OF THE BGG SIGNATURE}
\medskip

\def\mapright#1{\smash{ \mathop{\longrightarrow}\limits_{#1}}}
\def\maprightop#1{\smash{ \mathop{\longrightarrow}\limits^{#1}}}
\def\mapdown#1{\Big\downarrow\rlap{$\vcenter{\hbox{$\scriptstyle {#1}$}}$} }

Let ${\cal A_N}$ denote the set of all arrow-pairs of ${\frak S}_N$, as
defined above. In [BGG], it is proved that there exists a map
$$s:{\cal A_N}\to \{1,-1\}$$
with this property:
\smallskip
{\bf For every set of four arrow-pairs}
$$w_1\buildrel (i_1,j_1)\over \longrightarrow w_2,
w_2\buildrel (i_2,j_2)\over \longrightarrow w_4,
w_1\buildrel (i_3,j_3)\over \longrightarrow w_3,  
w_3\buildrel (i_4,j_4)\over \longrightarrow w_4,\eqno(3.1.7)$$
{\bf in} ${\cal A}$, {\bf we have}
$$s(w_1,w_2)s(w_2,w_4)s(w_1,w_3)s(w_3,w_4)=-1\eqno(3.1.8)$$
We shall call such a map a BGG-{\bf signature} for ${\cal A_N}$.
(Although the treatment in [BGG], and so some of the following 
discussion in this   section,
 make sense in much greater generality, for our present 
purposes it suffices to restrict to the consideration of ${\cal A_N}$. )
[BGG] calls a configuration such as (3.1.7) a 
{\bf square}; we denote (3.1.7) by the diagram
$$\matrix{
w_1&\mapright{(i_1,j_1)}&w_2\cr
\mapdown{(i_3,j_3)}&&\mapdown{(i_2,j_2)}\cr
w_3&\mapright{(i_4,j_4)}&w_4\cr
}\hskip.25in\eqno(3.1.9)$$ 
\smallskip
More generally, if ${\cal A'_N}$ is a subset of ${\cal A_N}$, by a 
{\bf partial} BGG-{\bf signature for} ${\cal A_N'}$, will be meant a map 
$$s:{\cal A'_N}\to \{1,-1\}$$
such that (3.1.8) holds for every square (3.1.9), for which all four 
arrow-pairs lie in ${\cal A_N'}$.
\medskip
The BGG resolutions, for the case $A_N$, are only completely specified
once one of the (many) BGG signatures for ${\cal A_N}$ has been selected.
Because Zelevinsky's construction in [Zel] 
of the complex studied in the present
section is obtained from the BGG resolution, it also is only specified up
to the choice of a BGG signature.
\medskip
This ambiguity is perhaps not very serious, but Akin in [Akin1,2] has shown
one specific way to choose these matters in a completely unique fashion
---let us pause here to sketch his normalization, since it is fairly 
short, and the author has not seen it
presented as an explicit algorithm in the literature.
Akin defines a specific BGG-signature, as follows:
\medskip
For every permutation $w\in {\frak S}_N$, let $P(w)$ 
denote the set of permutations $\overline w$ which precede $w$ 
in the Bruhat order---i.e., such that there exists a chain
$$w=w_0\buildrel (i_0,j_0)\over \longrightarrow w_1
\longrightarrow\cdots\longrightarrow w_l
\buildrel (i_l,j_l)\over\longrightarrow \overline w$$
of arrow-pairs in ${\frak S}_N$, beginning with $w$ and ending with 
$\overline w$. (We include $w$ in $P(w)$.)
\bigskip
Let
$$\tilde w={{1,2,\cdots,N-1,N}\choose{N,N-1,\cdots,2,1}}:i\mapsto N+1-i $$
denote the element in ${\frak S}_N$ of maximal length ${N\choose 2}$. 
 (Note that thus $P(\tilde w)={\frak A}_N$ and $P(I)
=\{I\}$.)       

Let $\Sigma$ denote a  chain
$$\Sigma:\tilde w=w_{{N \choose 2}}\to\dots\to w_0=I\eqno(3.1.10)$$
of ${N \choose 2}$ arrow-pairs in ${\frak S}_N$, such that each arrow-pair
$$w_{p+1}\buildrel (i_p,i_p+1)\over \longrightarrow w_p$$
in the chain $\Sigma$ is associated with an {\it elementary} transposition
 $(i_p,i_p+1)$. 
\smallskip
Then as Akin notes implicitly, the {\it existence
proof} in [BGG], actually gives a {\it constructive algorithm} for computing---
from the datum $\Sigma$---a BGG-signature $s_{\Sigma}$ for  ${\frak S}_N$, by
the inductive procedure next to be explained.
\medskip
Suppose that $0\le p<{N \choose 2}$, and that a partial BGG-signature  
$s_p$ has been given for $P(w_p)$. Then the argument in ([BGG],pp.56 and 57)
 in fact shows that the following rules furnish a well-defined extension of $s_p$
(involving no choices other than that of $\Sigma$)
 to a partial BGG-signature $s_{p+1}$ for $P(w_{p+1})$:
\medskip
Let $w\maprightop{(q,r)}w'$ be an arrow-pair, with 
$w$ (hence $w'$) in $P(w_{p+1})$.
In order to compute $s_{p+1}(w,w')$, we must consider four cases:

$\underline{\rm Case\, I:  w\in P(w_p)}$

Then also $w'\in P(w_p)$ Since we wish $s_{p+1}$ to 
extend $s_p$, we are forced to define
$$s_{p+1}(w,w')\buildrel \rm def \over=s_p(w,w')$$

$\underline{\rm Case\, II:  w\notin P(w_p),(q,r)=(i_p,i_p+1)}$ 

Here [BGG] defines $s_{p+1}(w,w')$ to be +1.(This {\it Ansatz}
implies the rules for the two following Cases.)

$\underline{\rm Case\, III:  w\notin P(w_p),w'\notin P(w_p), (q,r)\ne(i_p,i_p+1)}$

$\underline{\rm Case\, IV:  w\notin P(w_p),w'\in P(w_p), (q,r)\ne(i_p,i_p+1)}$

In both of these two cases, it follows from ([BGG], Lemma 11.3 on p.53), that
$$w\buildrel (i_p,i_p+1)\over \longrightarrow w(i_p,i_p+1)$$
and
$$w'\buildrel (i_p,i_p+1)\over \longrightarrow w'(i_p,i_p+1)$$ 
are arrow-pairs, and that
$w(i_p,i_p+1)$ and $w'(i_p,i_p+1)$  lie 
in $P(w_p)$.

Thus, in Case III, $s_{p+1}(w,w')$ is well-defined by
$$s_{p+1}(w,w')\buildrel \rm def \over
=-s_p(w(i_p,i_p+1),w'(i_p,i_p+1))$$     
(i.e., by 
the requirement, that the product of the edge-signatures of the following 
square
$$ \matrix{ w & \mapright{} & w'\cr \mapdown{+1}&&\mapdown{+1}\cr
w(i_p,i_p+1)  & \mapright{} & w'(i_p,i_p+1) \cr 
}$$ 
is to be $-1$.)

Similarly, in Case IV we must set
$$s_{p+1}(w,w')\buildrel \rm def \over
=-s_p(w(i_p,i_p+1),w'(i_p,i_p+1))s_p(w',w'(i_p,i_p+1)$$ 
\smallskip
Continuing inductively, this algorithm gives the desired BGG-signature
$s_{N(N-1)/2}=s_{\Sigma}$
on $P(\tilde w)={\cal A}_N$
\medskip
We shall (following Akin) utilize a canonical choice of the chain $\Sigma$,
 which is perhaps sufficiantly explained by giving the case $N=4$:
$$1234\longrightarrow\underbrace{1243}_{{\rm float}\, 3}\longrightarrow
\underbrace{1423\longrightarrow 1432}_{{\rm float}\, 2}\longrightarrow
\underbrace{4132\longrightarrow 4312 \longrightarrow 4321}_{{\rm float}\, 1}
\eqno (3.1.11)$$

The associated BGG-signature will be called the 
{\it Akin signature}, and will be denoted by $sgn_A$. 
\bigskip
\centerline{\bf C.\quad THE SHIFT-MATRICES SUBORDINATE TO AN ARROW-PAIR}

Let $(\pi,\pi')$ be an arrow-pair, so in particular there exist unique
integers $i,j$ satisfying (3.1.6). We define the {\it multiplicity} 
of $(\pi,\pi')$ to be the positive integer
$$r=\pi(i)-\pi(j)>0  \eqno(3.1.12)$$

We shall define $d_{\alpha}^{\pi,\pi'}$ to be a 
suitable ${\Bbb Z}$-linear combination of those
Weyl polarizations $P(\sigma)$, for which $\sigma$ lies in the set
defined as follows:

\proclaim Definition 3.1.1. Let 
$1\le i<j\le N$, and let r be a positive integer.
Then we denote by $TERM(i,j,r)$ the set of all $N$-shifts 
$$\sigma \in \Pi^N $$
which satisfy the following three conditions:
\item{I)}$\sigma_{p,q}=0$ unless $i\le q<p\le j$.
\item{II)}$r=\sum_{l=i+1}^j \sigma_{l,i}$
\item{III)}For all $k$ strictly between $i$ and $j$, we have
$$ \sum_{l=i}^j \sigma_{l,k} = \sum_{l=i}^j \sigma_{k,l} \eqno(3.1.13) $$
(i.e., $\sigma_{k,i}+\cdots+\sigma_{k,k-1}=\sigma_{k+1,k}+\cdots+
\sigma_{j,k}$.)\hfil \break
Also, it will be convenient to denote by $R_k(\sigma)$
the common value of both sides of eqn.(3.1.13).

In particular, if $i,j,r$ are related as above to an arrow-pair
 $(\pi,\pi')$ , then we shall say that the elements of the set
 $TERM(i,j,r)$ just defined, are {\it subordinate to} 
 $(\pi,\pi')$.

Let us again note specifically, that this condition on $\sigma$
 is totally independent of $\alpha$, and indeed only depends on $(i,j)$
 and $r$.
\smallskip
\centerline{{\bf D. THE NUMERICAL COEFFICIENT OF $P(\sigma)$
 IN } $d_l^{\pi,\pi'}$}
\smallskip
Let  $(\pi,\pi')$ be an 
arrow-pair, with $i,j$ determined by (2.2), and
with multiplicity $r=\pi(i)-\pi(j)$. Let 
$$\alpha=(a_1,...,a_N)\in {\Bbb C}^N\,,$$
and let $\sigma$ be an $N$-shift subordinate to  $(\pi,\pi')$.
 
\proclaim Definition 3.1.2. Under the preceding hypotheses, we define the
{\bf amplitude}
$$\langle \sigma;\pi,\pi' \rangle \in {\Bbb C}$$
to be the integer given by
$$\langle \sigma;\pi,\pi'\rangle=
r!\cdot \prod_{k=i+1}^{j-1} \left [
R_k(\sigma)!(r-R_k(\sigma) )!{\pi(i)-\pi(k) \choose r-R_k(\sigma)}
\right ] \eqno (3.1.14)$$              
(with the understanding the product in (3.1.14) is to be taken equal
to 1 if it is empty, i.e. if $j=i+1$.)

\medskip
Note:This amplitude is independent of $\alpha$, depending (as the notation
indicates) only on the N-shift $\sigma$ and the arrow-pair $(\pi,\pi')$.
This amplitude is an integer, which may be negative.
\bigskip
\centerline{\bf CONSTRUCTION OF THE MAPS $d_{\alpha}^{\pi,\pi'}$} 
We may now, finally, complete our construction of the natural transformations
$d_{\alpha}^{\pi,\pi'}$, and hence of the differentials in the complex
(3.1.5), as follows: 
\smallskip
Let  $(\pi,\pi')$ be an
arrow-pair, with $i,j$ determined by (2.2), and
with multiplicity $r=\pi(i)-\pi(j)$. Let
$$\alpha=(a_1,...,a_N)\in {\Bbb Z}^N\,,$$ 
let $V$ be a complex vectorspace, and let sgn be a BGG-signature
for ${\cal A}_N$(not necessarily the Akin signature). Finally, let
$$\omega\in ZEL^{\pi}(V,\alpha)=S^{b_1}V\otimes S^{b_2}V
\otimes \cdots \otimes S^{b_N}V$$
(with the $b$'s determined by (3.1.3a))
\medskip
Then we define
$$d_{\alpha,sgn}^{\pi,\pi'}\omega\buildrel \rm def \over{=}sgn(\pi,\pi')
\cdot\sum_{\sigma \in TERM(i,j,r)}
\langle \sigma;\pi,\pi' \rangle \cdot 
P(\sigma) \omega \eqno(3.1.15) $$     
\medskip
This completes our {\bf construction} of the complex
$ZEL(\alpha,V)$; we must postpone till Chapter 4 the {\bf proof}
of the following theorem, which asserts that the result of this
construction is indeed is an exact sequence, (except
possibly for its last term) and coincides completely
with the complex obtained by the method of Zelevinsky:

\proclaim Theorem 3.1.3. With 
$$d^{\pi,\pi'}_{\alpha}=d^{\pi,\pi'}_{\alpha,sgn}$$ 
given by (3.1.15),
 let us define the maps $d_k$ in (3.1.4) by
$$d_k=\bigoplus_{\scriptstyle l(\pi)=k \atop \scriptstyle l(\pi')=k-1
}d^{\pi,\pi'}_{\alpha}\;;$$ 
then these maps coincide completely with the maps
$$d_k^{ZEL}:\bigoplus_{\scriptstyle \pi \in {\frak S}_N
\atop l(\pi)=k}{\cal SYM}^{\alpha -\rho +\pi(\rho)}V 
 \to\bigoplus_{\scriptstyle \pi' \in {\frak S}_N
\atop l(\pi')=k-1}{\cal SYM}^{\alpha -\rho +\pi'(\rho)}V$$ 
 constructed by Zelevinsky
 in [Zel] from the BGG-resolutions for $A_N$ (constructed using sgn.)

Let us conclude this section by once again emphasizing, that the preceding
construction makes {\bf NO} use of the theory of Verma modules. By
contrast, the
{\bf proof} this construction yields an exact sequence---at 
least, the only complete proof available at
present to the author---makes heavy use of the work of
Verma, Bernstein-Gel'fand-Gel'fand,  Shapovalov and Zelevinsky.

We shall next consider some specific examples of the construction just
explained.

\def\az#1{ZEL[#1]}
\def\s#1#2#3{S^{#1}V\otimes S^{#2}V\otimes S^{#3}V}
\def\und#1{\underline{#1}}
\bigskip
\centerline{\S3.2 \bf SOME ILLUSTRATIVE EXAMPLES}
\bigskip
\centerline{\bf EXAMPLE 3.2.1}

Our first example involves the case N=3. The special case 
$N=3$ of   the Zelevinsky complex (with $\alpha$ a 
partition) is presented
in ([Doty], pp.134--136), where  this result is attributed to Verma. In Doty's
complex, the differentials are  explicitly furnished as elements of 
${\frak A}_3$, but not, of course, expressed in the language of
multi-polarizations. Thus, this earlier
data provides an excellent test for our assertions.
(Although it is assumed in [Doty]  that $\alpha$ is a partition, in
fact these results are valid without the assumption $a_1\ge a_2\ge a_3\ge 0$.)
As we shall see, the complex presented for $N=3$ by Doty and Verma, is
in fact precisely the Zelevinsky complex, normalized by the choice
of the Akin signature $sgn_A$.
\smallskip
Assume then,
$$N=3\,,\alpha=(a_1,a_2,a_3)\in {\Bbb Z}^3\;$$
and let V denote a   finite-dimensional complex vector-space.
It will be convenient to denote the permutation
$$\pi=\left ( {1 \atop \pi1}\,{2\atop \pi2}\,{3 \atop \pi3}\right)$$
by $[\pi1\, \pi2\, \pi3]$.
\smallskip
Here  the Zelevinsky complex $ZEL(\alpha,V)$ assumes the form
$$0\to ZEL^3\buildrel d_3\over \longrightarrow ZEL^2\buildrel d_2 \over
\longrightarrow ZEL^1 \buildrel d_1 \over \longrightarrow
ZEL^0\eqno(1)$$ 
with the terms $ZEL^l$ given by (3.1.3) and (3.1.4) as follows:
$$\eqalign{
ZEL^0&=S^{a_1}V\otimes S^{a_2}V \otimes S^{a_3} \cr
ZEL^1&=\az{213}\oplus\az{132} \cr
&=(\s{a_1+1}{a_2-1}{a_3})
\oplus (\s{a_1}{a_2+1}{a_3-1}) \cr
ZEL^2&=\az{231}\oplus\az{312} \cr
&=(\s{a_1+1}{a_2+1}{a_3-2})
\oplus(\s{a_1+2}{a_2-1}{a_3-1}) \cr
ZEL^3&=\az{321}=\s{a_1+2}{a_2}{a_3-2}\cr }$$

Having thus computed the {\it terms} in (1), let us 
next turn to the more interesting question of the {\it differentials}.
To obtain these, let us go in order through the four steps explained in 
\S3.1:

$\underline{\rm STEP\, A:}$ ${\cal A}_3$ consists of the 8 arrow-pairs:
$$\cases{
\tau_1:[213]\buildrel (1,2) \over \longrightarrow [123], 
\tau_2:[132]\buildrel (2,3) \over \longrightarrow [123], 
\tau_3:[231]\buildrel (1,3) \over \longrightarrow [132], 
\tau_4:[231]\buildrel (2,3) \over \longrightarrow [213],&\cr  
 \tau_5:[312]\buildrel (1,2) \over \longrightarrow [132],
\tau_6:[312]\buildrel (1,3) \over \longrightarrow [213],  
\tau_7:[321]\buildrel (1,2) \over \longrightarrow [231],  
\tau_8:[321]\buildrel (2,3) \over \longrightarrow [312]\;.&\cr
}\eqno(3.2.1)$$

$\underline{\rm STEP\, B:}$ We must next compute the Akin signature
$sgn_A$ for these 8 arrow-pairs. For chain (3.1.10) with $N$=3, we take
the Akin choice
$$w_3=[321]\buildrel (2,3) \over \longrightarrow 
w_2=[312]\buildrel (1,2) \over \longrightarrow w_1=    
[132]\buildrel (2,3) \over \longrightarrow w_0=I$$  
from which the inductive procedure of BGG yields the Akin
signatures given in the following table (where the arrow-pairs $\tau_i$
are given by (3.2.1)):
\smallskip
\centerline{
\vbox{\offinterlineskip\halign{
\hfil$#$\hfil&\vrule#\quad&#&\quad\vrule#\quad&#&\quad\vrule#\quad&#&
\quad\vrule#\quad&#&\quad\vrule#\quad&#&\quad\vrule#\quad&#&
\quad\vrule#\quad&#&\quad\vrule#\quad&#&\quad\vrule#\cr
\strut i&&1&&2&&3&&4&&5&&6&&7&&8&\cr  
\noalign{\hrule}
\strut sgn_A(\tau_i)&&+&&+&&$-$&&+&&+&&$-$&&+&&+&\cr
}}}
\smallskip
(For readers interested in applying the inductive  algorithm of \S3.1.1 to
verify this table, it may be helpful to note
that, setting
$$F_i=P(w_i)\backslash P(w_{i+1})\,,$$ 
there is induced on ${\frak S}_3$ the
filtration given by:
$$F_0=\{[123]\},F_1=\{[132]\},F_2=\{[213],[312]\},F_3=\{[231],[321]\}\;)$$
\smallskip
$\underline{\rm STEPS\,C
\, and\, D:}$ Here the maps 
$$d^{\pi,\pi'}_{\alpha}=d_{\alpha}(\tau)$$
(where $\tau$ denotes the arrow-pair $\pi \to \pi'$) which are
furnished by the remaining two steps, are given in 
the second row of the following table, while the third row lists
the maps $d_{\rm DV}(\tau)$ furnished by Doty and Verma; it is 
asserted that the second and third rows coincide, i.e., that
$d_{\alpha}(\tau)=d_{\rm DV}(\tau)$.
\smallskip
\centerline{
\vbox{\offinterlineskip\halign{
\hfil$#$\hfil&\vrule#&\hfil#\hfil&\vrule#&\hfil#\hfil&\vrule#
&\hfil#\hfil&
\vrule#&\hfil#\hfil &\vrule#&\hfil#\hfil&\vrule#&\hfil#\hfil&\vrule#\cr
\strut i&&1,8&&2,7&&3&&4&&5&&6&\cr
\noalign{\hrule}
\strut d_{\alpha}(\tau_i)&&$E_{2,1}$&&$E_{3,2}$&&
$E_{3,1}-P(E_{3,2}+E_{2,1})$&&$2P(2E_{3,2})$&&$2P(2E_{2,1})$&&
$-2E_{3,1}-P(E_{3,2}+E_{2,1})$&\cr
\noalign{\hrule}
\strut d_{\rm DV}(\tau_i)&&$E_{2,1}$&&$E_{3,2}$&&$E_{3,2}E_{2,1}-
2E_{2,1}E_{3,2}$&&$(E_{3,2})^2$&&$(E_{2,1})^2$&
&$E_{2,1}E_{3,2}-2E_{3,2}E_{2,1}$&\cr
}}}   
\smallskip

We chase through
the details for the case
$$\tau_3:[231]\buildrel (1,3) \over \longrightarrow [132]\,;$$
it is left as an exercise to 
any reader so interested, to verify the other 
entries in this table by the same straightforward algorithm.

The natural transformation
$$\eqalignno{
&d_{\alpha}(\tau_3):ZEL[2,3,1]=S^{a_1+1}V\otimes 
S^{a_2+1}V\otimes S^{a_3-2}V&(3.2.2)  \cr 
&\to ZEL[1,3,2]=S^{a_1}V\otimes S^{a_2+1}V\otimes S^{a_3-1}V&\cr }$$
is one of the four constituents of $d_2$. 

Here 
$$\pi={123 \choose 231},\pi'={123 \choose 132},(i,j)=(1,3)$$
Hence by (3.1.8), $r=\pi(1)-\pi(3)=1$. Examination of Def.3.1.1 
(in Case C of the preceding subsection)
shows that the 3-shift $\sigma$ is subordinate to $\pi$
 and $\pi'$, (i.e., lies in
$TERM(1,3,1)$), precisely when:
\item{(i)} All $\sigma_{i,j}$ are 0, except possibly $\sigma_{2,1}, 
\sigma_{3,1}, \hbox{\rm and}\; \sigma_{3,2}$.
\item{(ii)} $\sigma_{2,1}+\sigma_{3,1}=1$
\item{(iii)} $\sigma_{2,1}=\sigma_{3,2}$

Hence
$$TERM((1,3,1)=\{\sigma_1,\sigma_2\}$$
with 
$$\sigma_1=E_{3,1}\,,\sigma_2=E_{2,1}+E_{3,2}$$ 
Thus, $d_{\alpha}(\tau_3)$ 
is a ${\Bbb C}$-linear combination of $P(2E_{3,1})$ and
$P(E_{2,1}+E_{3,2} )$, with the numerical coefficients next to
be determined (via Def.3.1.2 in Step D of the preceding subsection.) 

Consider first $$\sigma_1=E_{3,1}\;.$$
Here $$R_2(\sigma_1)=0\;,$$ and eqn.(3.1.15) yields the amplitude
$$<E_{3,1};\pi,\pi'>=1!0!1!{\pi(1)-\pi(2) \choose 1}={-1 \choose 1}=-1$$

Similarly, we have $$R_2(\sigma_2)=1\;,$$and so
$$<E_{3,2}+E_{2,1};\pi,\pi'>=1$$
Recalling that $sgn_A(\tau_3)=-1$, (3.1.16) gives the entry
$$d_{\alpha}(\tau_3)=E_{3,1}-P(E_{3,2}+E_{2,1})\eqno(3.2.3a)$$
in the preceding table.

It remains to verify that this map coincides, on their common domain
$$ZEL[2,3,1]=S^{a_1}V\otimes S^{a_2+1}V\otimes S^{a_3-1}V\;,$$
with the Doty-Verma map
$$d_{\rm DV}(\tau_3)=[E_{3,2}E_{2,1}-2E_{2,1}E_{3,2}]|ZEL[2,3,1]\eqno(3.2.3b)$$
To prove this, it suffices to verify that both maps have the same
effect on every element
$$\eqalignno{
\omega&=(x_1\cdots x_{a_1+1})\otimes(y_1\cdots y_{a_2+1})\otimes
(z_1\cdots z_{a_3-2})\cr
\noalign{\hbox{in}}
ZEL[231]&=S^{a_1+1}V\otimes 
 S^{a_2+1}V\otimes S^{a_3-2}V\cr } $$
(all $x$'s, $y$'s and $z$'s in $V$). Let us set
$$
A:=E_{3,1}\omega=\sum_{i=1}^{a_1+1}(x_1\cdots\widehat{x_i}\cdots x_{a_1+1})
\otimes(y_1\cdots y_{a_2+1})\otimes(x_i\cdot z_1\cdots z_{a_3-2})$$
and
$$\eqalign{
&B:=P(E_{3,2}+E_{2,1})\omega= \cr
&=\sum_{i=1}^{a_1+1}\sum_{j=1}^{a_2+1}
(x_1\cdots\widehat{x_i}\cdots x_{a_1+1})\otimes
(x_i\cdot y_1\cdots\widehat{y_j}\cdots y_{a_2+1})\otimes
(y_j\cdot z_1\cdots z_{a_3-2})\cr}
$$
(Note that $A$ and $B$ both lie in 
$ZEL[1,3,2]=S^{a_1}V\otimes S^{a_2+1}V\otimes S^{a_3-1}V$.)

Then direct computation (as explained in \S1.2) shows that
$$E_{3,2}E_{2,1}\omega=E_{3,2}\sum_{i=1}^{a_1}(x_1\cdots\widehat
x_i\cdots x_{a_1+1})\otimes(x_i\cdot y_1\cdots y_{a_2+1})\otimes
(z_1\cdots z_{a_3-2})=A+B$$
and similarly
$$E_{2,1}E_{3,2}\omega=B$$
Hence, finally,
$$d_{\rm DV}(\tau_3)\omega=(A+B)-2B=A-B=E_{3,1}\omega-P(E_{3,2}+E_{2,1})\omega
=d_{\alpha}(\tau_3)\omega\;.$$
\medskip
($\underline{\rm Question}$: Is there some simple rule one could use
directly to compute the Akin signature---or some other specific
BGG-signature---rather than the tedious step-by-step inductive algorithm
presented in \S3.1?)

\centerline{\bf EXAMPLE 3.2.1} 
Our next example involves $N=4$ and the arrow-pair 
$$\tau:[2341]\buildrel (2,4)\over \longrightarrow [2143]$$
Let $sgn$ denote an arbitrary choice (not necessarily $sgn_A$) of 
BGG-signature on ${\cal A}_4$.

If $\alpha=(a_1,a_2,a_3,a_4)\in{\Bbb Z}^4$, and $V$ varies over complex
vector-spaces, then the natural transformation $d_{\alpha,V}(\tau)$, which
maps
$$ZEL^{[2341]}(\alpha,V)=S^{a_1+1}V\otimes S^{a_2+1}V\otimes 
S^{a_3+1}V\otimes S^{a_4-3}V$$
into
$$ZEL^{[2143]}(\alpha,V)=S^{a_1+1}V\otimes S^{a_2-1}V\otimes S^{a_3+1}V
\otimes S^{a_4-1}V\;,$$
is a component of the differential $d_3$ in $ZEL(\alpha,V)$ (since
$[2341]$ has 3 inversions.) This mapping is, in fact, that resulting from the
action on the ${\frak A}_4$-module $ZEL^{[2341]}(\alpha,V)$ of the following 
element in ${\frak A}_4$ (which, it should be noted, is completely
independent both of $\alpha$ and of $V$):
$$d(\tau)=sgn(\tau)\cdot[4P(2E_{3,2}+2E_{4,3})-2P(E_{3,2}+E_{4,3}+E_{4,2})
+2P(2E_{4,2})]\eqno(3.2.5)$$
(The reader may find it a helpful exercise, to compute
(3.2.5), using the algorithm explained in \S3.1.)
\medskip
One final bit of propaganda for the efficacy and appropriateness of
the Weyl polarizations in such computations: let us examine the precise effect 
of (3.2.5) on the generating element
$$\omega=\und{w}\otimes\und{x}\otimes\und{y}\otimes\und{z}=
(w_1\cdots w_{a_1+1})\otimes(x_1\cdots x_{a_2+1})\otimes
(y_1\cdots y_{a_3+1})\otimes(z_1\cdots z_{a_4-3})
$$
for $ZEL^{[2341]}(\alpha,V)$:

Using eqn.(3.2.4) we obtain
$$d_{\alpha,V}(\tau)\omega=sgn(\tau)(4A-2B+2C)$$
where we have set:
$$A=P(2E_{3,2}+2E_{4,3})\omega=
\sum_{\scriptstyle {i<i'}\atop{i,i'\in\und{a_2+1}}}
\sum_{\scriptstyle {j<j'}\atop{j,j'\in\und{a_3+1}}}
\und{w}\otimes {\und{x}\over{x_i\cdot x_{i'}}}\otimes
(x_i\cdot x_{i'}\cdot{\und{y}\over{y_j\cdot y_{j'}}})\otimes 
(y_j\cdot y_{j'}\cdot\und{z})\;,$$
$$B=P(E_{3,2}+E_{4,3}+E_{4,2})\omega=
\sum_{\scriptstyle {i\ne i'}\atop{i,i'\in\und{a_2+1}}}
\sum_{j\in\und{a_3+1}}\und{w}\otimes{\und{x}\over{x_i\cdot x_{i'}}}
\otimes(x_i\cdot{\und{y}\over y_j})\otimes(x_{i'}\cdot y_j\cdot\und{z})$$
and
$$C=P(2E_{4,2})\omega=\sum_{\scriptstyle {i<i'}\atop{i,i'\in\und{a_2+1}}}
\und{w}\otimes {\und{x}\over{x_i\cdot x_{i'}}}\otimes
\und{y}\otimes(x_i\cdot x_{i'}\cdot\und{z})$$
Note that $A,B,C$ all lie, as they ought to, in
$$ZEL^{[2143]}(\alpha,V)=S^{a_1+1}V\otimes S^{a_2-1}V\otimes S^{a_3+1}V
\otimes S^{a_4-1}V\;,$$

\bigskip
\centerline{\S3.3 \bf The action on $S^{\otimes N}$ of the
Verma-Shapovalov Elements for $A_N$} 
\medskip
Let ${\frak g}$ be a semi-simple complex Lie algebra, with
${\frak h}$ a selected Cartan subalgebra; also assume selected an ordering for
($\frak g,\frak h$), with $\Delta^+$ the set of positive roots. Let
$\rho$ denote half the sum of the
roots in $\Delta^+$. For any positive root
$\alpha$, we denote the associated co-root in $\frak h$ by $h_{\alpha}$
 ---so that 
$$s_{\alpha}(\lambda)=\lambda-\lambda(h_{\alpha})\alpha$$ 
 for all $\lambda\in {\frak h}^*$. Let
$\frak N_+$ denote the nilpotent sub-algebra of $\frak g$ generated by the
positive root spaces, and $\frak N_-$ that generated by the negative 
root spaces.
\par\noindent
\medskip
For any $\lambda \in \frak h^{*}$, we denote by ${\cal I}_{\lambda}$,
the left ideal in the enveloping algebra $\frak A$($\frak g$) of 
 $\frak g $, generated by
$$\Delta^+ \cup \{h-\lambda(h)\cdot 1 : h \in \frak h \}$$
Thus the $\frak g$-module
$${\cal V}_{\lambda}:=\frak A(\frak g)/{\cal I}_{\lambda}$$
is precisely the Verma module over $\frak g$, with highest weight
$\lambda$. We denote by $v_\lambda$ the image of 1 in this quotient,
(so that ${\cal V}_{\lambda}$ is cyclic on the distinguished highest-weight
vector $v_{\lambda}$.)
\medskip
\proclaim Definition 3.3.1.  By a {\bf Verma triple} for such a
${\frak g}$, will 
be meant an ordered triple ($\alpha,r,\lambda$) satisfying the
following four conditions:
\smallskip
\item{(i)} $\alpha$ is a positive root in ($\frak g,\frak h$).
\item{(ii)} $\lambda \in \frak h^{*}$
\item{(iii)}$r$ is a positive integer
\item{(iv)}$\lambda-s_{\alpha}\bullet \lambda =r\alpha$
\par
{\bf REMARK}: Here the symbol $\bullet$ designates the ``twisted
action'' of the Weyl group on $\frak h$, given by
$$s_{\alpha}\bullet \lambda=s_{\alpha}(\lambda+\rho)-\rho$$
Thus condition (iv) may be replaced by the equivalent condition
$$(\lambda + \rho)h_{\alpha}=r$$


\bigskip
It is a well-known result, due to 
Verma, that given such a triple,
there exists a non-zero $\frak g$-linear homomorphism
$${\cal V}_{\lambda - r\alpha} \to {\cal V}_{\lambda}\eqno(3.3.1)$$
unique up to scalar multiples. However, since our purpose here is
to express this map as an explicit 
${\Bbb Z}$-linear combination of (the actions of)
Weyl polarizations, it is necessary
to select an explicit {\bf normalization} of this element; Shapovalov
has constructed one method for doing so, as follows:

\noindent
\proclaim Definition 3.3.2. Let $\tau$=($\alpha,r,\lambda$) be a Verma triple
for ${\frak g}$. \hfil \break 
 By the \hbox{{\bf Verma-Shapovalov element}}
$$\gamma=VS(\tau)=VS_{\alpha,r}(\lambda)$$
for this triple, is meant the element $\gamma$ in 
$\frak A(\frak N_- )$ uniquely specified by  the two following conditions:
\smallskip
VS1)There is a non-zero $\frak g$-linear map
$$\phi:{\cal V}_{\lambda - r\alpha} \to {\cal V}_{\lambda}$$
such that
$$\phi(v_{\lambda - r\alpha})=\gamma v_{\lambda}$$
\smallskip
VS2)Choose a total ordering $<<$ of $\Delta_+$; say
$$\{\alpha_1<<\alpha_2<<\cdots<<\alpha_m\}\hbox{,  where }m=\#(\Delta^+)$$
For each $\alpha\in\Delta_+$ let $E_{\alpha}$ be an associated root vector.
\medskip
This choice associates to every map
$\pi:\Delta^+ \to {\Bbb N}$
a basis element
$$F_{\pi}:=(E_{\alpha_1})^{\pi(\alpha_1)}\cdots (E_{\alpha_m})^{\pi(\alpha_m)}
\eqno(3.3.2)$$
in the PBW basis for $\frak A(\frak N_-)$ associated with $<<\,$, 
 namely the basis
$$\{F_{\Pi}|\Pi \in (\Bbb N)^{\Delta^+} \}\,. \eqno(3.3.3) $$
Note in particular the distinguished map
$$\pi<r>:\Delta^+ \to {\Bbb N}$$
which maps each of the simple roots
to r, and maps all other positive roots to 0.
It is then required that,  when the Verma-Shapovalov element 
$\gamma$ is expanded
as a $\Bbb C$-linear combination of the PBW-basis (3.3.3),
the basis vector $F_{\pi<r>}$ shall have coefficient precisely~1.

\medskip 
\noindent{\bf REMARKS:}\quad The first requirement 
VS1), specifies $\sigma=VS(\tau)$
uniquely, up to a scalar multiple. (This is the fundamental result of Verma
 which is the basis of the present paper.) Concerning the second
requirement, (due it seems to Shapovalov),
whose effect is to remove this last ambiguity in 
the definition of VS($\tau$), let us note that this requirement 
presupposes two non-obvious facts:
1)The coefficient of $F_{\pi[\tau]}$ is not identically 0 in all 
operators satisfying VS1). 
2)The requirement VS2) is in fact independent of 
the particular choice of total ordering $<<$ on $\Delta^+$.

These two facts may be found proved, in Franklin([Fra,\S3 and 4]
\footnote*{\sevenrm Caution: There is a typo in the statement of VS2) on
p.66 of [Fra]; in (3),loc.cit., $r=\sum n_i{\epsilon}_i$ must be
replaced by  $dr=\sum n_i{\epsilon}_i$ } 
for arbitrary semi-simple Lie algebras, in
arbitrary characteristic.

In a recent clarifying discussion of these matters, the present author 
was informed by Verma,
 that the construction in Verma's thesis [Verma 1] of a
${\frak g}$-homomorphism (3.3.1) was
in fact defined
 absolutely, not simply up to scalar multiples. This would imply
that the construction in Verma's thesis supplied an 
intrinsic method of normalizing (3.3.1).
Below, there will be given a formula (3.3.6), an expression in
which will then be verified to 
satisfy the two Shapovalov conditions listed above, and hence to coincide
with the Verma-Shapovalov element.The present author does 
not know the precise relation
between these two methods of normalization (Verma's
 and Shapovalov's), and hence must here leave open the question
of the relations (if any) of formula (3.3.6) below, to 
the normalization of (3.3.1)
propounded by Verma.

\bigskip
{\bf From now on, for the remainder of the present paper,
we shall restrict ourselves entirely to the Lie algebra}
$$\frak g = {\frak sl}_N(\Bbb C) =A_{N-1}  $$
{\bf and to its enveloping algebra} $(\frak A_N )^0$. 

In this special case, we choose, (as is usual), the  Cartan subalgebra
$\frak h$, to  be that formed by the diagonal $N\times N$ matrices of trace 0.
The dual  $\frak h^{*}$ of this is spanned over $\Bbb C$ by the $N$
linear functionals $\lambda_i$ ( $i$ between 1 and $N$) where
$\lambda_i$ maps any $N\times N$ matrix C in $\frak g$ into its
$i$-th diagonal entry $C_{i,i}$. (Thus the $\lambda_i$ have sum 0. )
The usual choice here for the set of positive roots is
$$\Delta^+=\{\lambda_i-\lambda_j : 1\le i< j \le N \}$$
and the associated nilpotent subalgebra $\frak N_+$ is that formed by the
$N\times N$ strictly upper-triangular matrices over $\Bbb C$---so
that  $\frak N_-$ is the Lie algebra formed by the strictly
lower-triangular ones.

With these choices, it is readily seen that the Verma triples
for $A_{N-1}$ consist of those ordered triples
$$(\alpha=\lambda_i-\lambda_j,r,\lambda=\sum_{i=1}^N l_i \lambda_i)$$
(with $1\le i<j\le N,r \in \Bbb Z^+ \hbox{, the } l_i$ being complex
numbers), which satisfy
$$l_i-l_j-i+j=r \eqno(3.3.4)$$

\proclaim Definition 3.3.3. Let
$$\tau=(\lambda_i-\lambda_j,r,\lambda=\sum_{i=1}^N l_i {\lambda}_i)$$
be a Verma triple for $A_{N-1}$, and (using Def.3.1.1) let 
$$\sigma \in TERM(i,j,r) \,;$$
Then we define the {\bf amplitude}
$\langle \sigma;\tau \rangle \in {\Bbb C}$ to be
$$\langle \sigma;\tau \rangle = r!\cdot \prod_{k=i+1}^{j-1} \left [
R_k(\sigma)!(r-R_k(\sigma))!{l_i-l_k \choose r-R_k(\sigma)}
\right ] \eqno (3.3.5)$$ \par

Our goal in the next \S\ will be the proof of:
\proclaim Theorem 3.3.4. If   
$$\tau=(\lambda_i-\lambda_j=\alpha,r,\lambda)$$
is a Verma triple for $A_{N-1}$, and if we set
$${\cal} T(\tau)\buildrel\rm def\over = 
TERM(i,j,r)$$
then 
$$VS(\tau)=\sum_{\sigma \in {\cal T}({\tau})}\langle \sigma;\tau \rangle
P(\sigma)\,.  \eqno (3.3.6) $$
\par
\bigskip
\noindent {\bf NOTE:} As noted by K.Akin in [Akin2, p.418], the maps 
$$d_{\alpha,sgn}^{\pi,\pi'}$$ for the Zelevinsky complex, are 
given by precisely the same elements
$$sgn(\pi,\pi')VS(\pi,\pi')$$
in ${\frak A}_N$ which furnish the corresponding maps in the BGG-resolution.
Hence Th.3.3.4 implies (and is apparently rather stronger than) Th.3.1.3.
\bigskip
\centerline{\bf Section 4\qquad Proof of the Assertions in Section 3}
\medskip
As observed at the end of \S3.3, 
to prove all the  assertions in \S3, it suffices to prove Th.3.3.4.
\smallskip
For the rest of this \S, there will be assumed
 the hypotheses of Th. 3.3.4 ; that
is, we suppose:
\item{(i)} $\tau=(\lambda_i -\lambda_j,r,\lambda=\sum l_i \lambda_i)$ 
is a Verma triple
for $\frak A_{N-1}$;\ and we set
\item{(ii)}${\cal T}(\tau)=TERM(i,j,r)\,,$
as defined in Def.3.1.1.
\smallskip
Note that (i) means that:\hfil \break
$1\le i<j\le N$, r is a positive integer, and
$$l_i-l_j-i+j=r \eqno(4.1)$$
\smallskip
Then, to complete the proof of Theorem 3.3.4, it suffices to prove
that the element $\gamma$ in $\frak A_N$, defined by
$$\gamma\buildrel\rm def\over =\sum_{\sigma \in {\cal T}({\tau})}
\langle \sigma;\tau \rangle P(\sigma) \eqno(4.2) $$ 
satisfies (relative to the given Verma triple $\tau$) the two conditions
VS1) and VS2) which, in Def.3.2.2, characterize the Verma-Shapovalov element
$VS(\tau)$.(Here the coefficients
$$\langle \sigma;\tau \rangle$$ 
 are the numbers given by
Def.3.3.3.)

Unwrapping all this, we see that the goal of proving Th.3.3.4
will have been achieved, once the three following assertions have
been established:
\item{A1}$\gamma\cdot v_{\lambda}$ is a maximal vector in 
${\cal V}_{\lambda}$.
\item{A2}This vector $\gamma\cdot v_{\lambda}$ has weight
$\lambda -r\alpha $.
\item{A3} $\gamma$ satisfies the normalization condition $VS2)$ in Def.3.3.2.

(Note that, if A1 and A2 are satisfied, then the element
$\gamma\cdot v_{\lambda}$ in ${\cal V}_{\lambda}$ 
generates a sub-module isomorphic to ${\cal V}_{\lambda-r\alpha}$,
from which VS1) is immediate.)

The next three sub-sections are devoted to the proof, in
this order, of these three assertions.

\def\ignore#1{}
\bigskip
\centerline{\bf \S4.1  Proof That $\gamma\cdot v_{\lambda}$ is a
Maximal Vector in ${\cal V}_{\lambda}$}
\medskip
The purpose of this sub-section is the proof of assertion A1, that 
is, the proof that:
$$E_{p,p+1}\cdot\gamma\cdot v_{\lambda}=0 \hbox{ for }1\le p\le N-1
 \eqno (4.1.1)$$
            
In the remainder of this paper, $ \equiv $ is always to be understood to mean
$\equiv$ modulo the left ideal ${\cal I}_{\lambda}$ in
${\frak A}_N$ defined 
in $\S 3.1$.
In other words, for 
all $\gamma$ and $\gamma'$ in ${\frak A}_N$,
$$\gamma\equiv\gamma'\hbox{  if and only if  }\gamma\cdot v_{\lambda}= 
\gamma'\cdot v_{\lambda}\,.$$
For instance, eq.(4.1.1) is equivalent to the assertion
$$E_{p,p+1}\cdot \gamma\equiv 0 \eqno(4.1.2)$$
This in turn is equivalent to the assertion  that
$$[E_{p,p+1},\gamma] \equiv 0 \hbox{ for }1\le p\le N-1 \eqno(4.1.3)\,, $$
since by definition every left multiple in ${\frak A}_N$ of
$E_{p,p+1}$ lies in  ${\cal I}_{\lambda}$.

\medskip

\centerline{$\bullet$ Throughout the rest of this sub-section, $1\le p<N$.}
\bigskip

\medskip
\noindent Our proof of (4.1.3) will proceed in three steps: \hfil \break
We first analyze the commutators
$$[E_{p,p+1},P(\sigma)]$$
for all $\sigma$ in $TERM(\tau)$. The results thus obtained             
will be applied in the second
step to expand $[E_{p,p+1},\gamma]$ as
a linear combination of Weyl polarizations with explicitly described
integer coefficients. In
the third step, we shall finally prove (4.1.3) by showing these coefficients
are all 0.

\medskip
\centerline{\bf Step One: Analysis of $[E_{p,p+1},P(\sigma)]$ for
$\sigma\in TERM(\tau)$}
 
\noindent By  Cor. 2.2.3, we have
$$[E_{p,p+1},P(\sigma)]=A-B \eqno(4.1.4)$$
where
$$A=\sum_{k=1}^N(\sigma_{p,k}+1)P(\sigma +E_{p,k}-E_{p+1,k})
\eqno(4.1.4a)$$
and
$$B=\sum_{k=1}^N (\sigma_{k,p+1}+1)P(\sigma +E_{k,p+1}-E_{k,p})
\eqno(4.1.4b) $$
\bigskip
By hypothesis, $\sigma$ satisfies the conditions of 
Def. 3.1.1. (in Part C of \S3.1).
In particular, condition I) of this definition implies that
$$\sigma_{p+1,k}=0 \hbox{ if }i\le k<p+1\le j \hbox{ does not hold,}$$
in which case 
$$(\sigma+E_{p,k}-E_{p+1,k})_{p+1,k}=-1\,,$$
so by eq.(1.4.8),
$$P(\sigma+E_{p,k}-E_{p+1,k})=0\,.$$
Since also $\sigma_{p,p}=0$, we obtain
$$A=1\cdot P(\sigma-E_{p+1,p}+E_{p,p})+
\sum_{k=i}^{p-1}(\sigma_{p,k}+1)(P(\sigma-E_{p+1,k}+E_{p,k})$$
In the same way, we have
$$B=P(\sigma-E_{p+1,p}+E_{p+1,p+1})+
\sum_{k=p+2}^j(\sigma_{k,p+1}+1)P(\sigma-E_{k,p}+E_{k,p+1})$$ 
Combining the two preceding equations with (4.1.4), we obtain:
$$\eqalignno{
[E_{p,p+1},P(\sigma)]&=[P(\sigma-E_{p+1,p}+E_{p,p})
-P(\sigma-E_{p+1,p}+E_{p+1,p+1})]+  &(4.1.5) \cr
&+\sum_{k=i}^{p-1}(\sigma_{p,t}+1)(P(\sigma-E_{p+1,k}+E_{p,k})- \cr
&-\sum_{k=p+2}^j(\sigma_{k,p+1}+1)P(\sigma-E_{k,p}+E_{k,p+1}) \cr
}$$
\medskip
\noindent It is next claimed that (4.1.5) is 0 unless $i\le p<j$:\hfil\break

\noindent Indeed, if $p\le i-1$ then both $\sigma-E_{p+1,p}+E_{p,p}$
 and $\sigma-E_{p+1,p}+E_{p+1,p+1}$ are non-effective (since
 $\sigma_{p+1,p}=0$); moreover, the first sum $\sum_{k=i}^{p-1}$
 in (4.1.5) is empty (hence 0), while in the second sum, each term is 
 0 (since each $\sigma-E_{k,p}+E_{k,p+1}$ is non-effective)---hence
 (4.1.5) is 0 in this case. Similarly, (4.1.5) is 0 if $p\ge j$.
\smallskip
\noindent We may thus (without loss of generality) strengthen as follows our 
earlier assumption:

\centerline{$\bullet$ Throughout the rest of this sub-section, $i\le p<j$.}
\bigskip
For all $r$ between 1 and $N$, let us set
$$col_r(\sigma)=(r^{th}\hbox{ column-sum of }\sigma)=\sigma_{1,r}+
\sigma_{2,r}+\cdots+\sigma_{N,r}$$
Because of our hypothesis that 
$$\sigma \in TERM(\tau)\,,$$
we have 
$$col_r(\sigma)=\cases{d,&if $r=i$;\cr
R_r(\sigma),&if $i<r<j$;\cr
0&otherwise.\cr}\eqno(4.1.6)
$$
\medskip
Consider next the first expression
$$F=[P(\sigma-E_{p+1,p}+E_{p,p})
-P(\sigma-E_{p+1,p}+E_{p+1,p+1})] \eqno(4.1.7)$$
occurring on the right side of eq.(4.1.5).
Since all diagonal entries of $\sigma$ vanish (because of condition
I in Def. 3.2.1), it then follows from Prop.2.2.2 that
$$P(\sigma-E_{p+1,p}+E_{p,p})=P(\sigma-E_{p+1,p})E_{p,p}
-(col_p(\sigma)-1)P(\sigma-E_{p+1,p})$$
 (if $\sigma-E_{p+1,p}$ is effective---but this equation is also valid
if $\sigma-E_{p+1,p}$ is not effective, since then both sides are 0)

Similarly, Prop.2.2.2 implies that
$$P(\sigma-E_{p+1,p}+E_{p+1,p+1})=P(\sigma-E_{p+1,p})E_{p+1,p+1}
-(col_{p+1}(\sigma))P(\sigma-E_{p+1,p})$$
Combining these two equations, with the fact that the left ideal
${\cal I}_{\lambda}$ contains (for all $k$ between 1 and $N$)the elements
$$E_{k,k}-E_{k+1,k+1}-\lambda(E_{k,k}-E_{k+1,k+1})=
E_{k,k}-E_{k+1,k+1}-l_k+l_{k+1}$$
we obtain the following congruence for (4.1.7): modulo ${\cal I}_{\lambda}$,
$$F \equiv (l_p-l_{p+1}-
col_p(\sigma)+col_{p+1}(\sigma)+1)P(\sigma-E_{p+1,p})$$

Inserting this last congruence into (4.1.5),we obtain
$$\eqalignno{
[E_{p,p+1},P(\sigma)]&\equiv(l_p-l_{p+1}-
col_p(\sigma)+col_{p+1}(\sigma)+1)P(\sigma-E_{p+1,p}) &(4.1.8) \cr
&+\sum_{k=i}^{p-1}(\sigma_{p,k}+1)(P(\sigma-E_{p+1,k}+E_{p,k})- \cr
&-\sum_{k=p+2}^j(\sigma_{k,p+1}+1)P(\sigma-E_{k,p}+E_{k,p+1}) \cr
}$$
\bigskip
 
\centerline{\bf Step Two: Expansion of $[E_{p,p+1},\gamma]$ into 
Multi-polarizations}
\bigskip
We combine (4.1.8), with
$$[E_{p,p+1},\gamma]=\sum_{\sigma\in TERM(\tau)}
\langle \sigma;\tau \rangle [E_{p,p+1},P(\sigma)]$$
to obtain the congruence 
$$[E_{p,p+1},\gamma]\equiv{\cal S}'_p+{\cal S}''_p+{\cal S}'''_p\,, 
\eqno(4.1.9) $$
where
$$\eqalignno{
{\cal S}'_p&=\sum_{\sigma\in TERM(\tau)}
\langle \sigma;\tau \rangle (l_p-l_{p+1}-col_p(\sigma)+
col_{p+1}(\sigma)+1)P(\sigma-E_{p+1,p}) &(4.1.9a) \cr
{\cal S}''_p&=\sum_{\sigma\in TERM(\tau)}\sum_{k=i}^{p-1}
\langle \sigma;\tau \rangle
(\sigma_{p,k}+1)P(\sigma-E_{p+1,k}+E_{p,k}) &(4.1.9b) \cr 
&&{\rm and}\cr
{\cal S}'''_p&=-\sum_{\sigma\in TERM(\tau)}\sum_{k=p+2}^j
\langle \sigma;\tau \rangle
(\sigma_{k,p+1}+1)P(\sigma-E_{k,p}+E_{k,p+1}) &(4.1.9c) \cr}
$$ 

\medskip
We must next collect coefficients of the various P's in (4.1.9):
\proclaim Lemma 4.1.1. 
$$[E_{p,p+1},\gamma]\equiv\sum_{\sigma\in TERM(\tau)}
\{A_p(\sigma)+B_p(\sigma)+C_p(\sigma)\}
P(\sigma-E_{p+1,p}) \eqno(4.1.10)$$
where
$$\eqalignno{
A_p(\sigma)&=\sum_{\sigma\in TERM(\tau)}((l_p-l_{p+1}-col_p(\sigma)+
col_{p+1}(\sigma)+1)\langle \sigma;\tau \rangle &(4.1.10a)\cr
B_p(\sigma)&=\sum_{k=i}^{p-1}\sigma_{p,k}
\langle \sigma+E_{p+1,k}-E_{p,k}-E_{p+1,p};\tau \rangle &(4.1.10b)\cr
&&{\rm and} \cr
C_p(\sigma)&=-\sum_{k=p+2}^j\sigma_{k,p+1}
\langle \sigma+E_{k,p}-E_{k,p+1}-E_{p+1,p};\tau \rangle &(4.1.10c) \cr
}$$

{\bf Note:} As explained earlier, in the present paper we adopt
everywhere the two conventions, that
$$\sigma\not\in  TERM(\tau)\Rightarrow \langle\sigma;\tau
\rangle=0\,.$$
and that
$$\sigma\in\Pi_{\pm}^N \hbox{ and $\sigma$ non-effective, imply  }
P(\sigma)=0$$
In particular, the preceding equations are to be interpreted in this way.
Thus, in the ensuing argument, some caution will be needed, to 
distinguish
the `regular' terms in these sums, from those which are assigned the value 
0 by the conventions just reviewed.

\smallskip
{\bf Proof of Lemma 4.1.1:}
To prove the Lemma, it clearly suffices to verify the following two assertions:

$\underline{{\bf CLAIM\, ONE:}}$ {\it If} $\sigma\in TERM(\tau)$
{\it then } $P(\sigma-E_{p+1,p})$
 {\it  either is 0, or has the same coefficient
 on both sides of (4.1.10).}

$\underline{\bf CLAIM \, TWO:}$ {\it If $P(\sigma')$ occurs with
nonzero coefficient as a term  in at least one of the three sums
 (4.1.10a, b, c), then} $\sigma'+E_{p+1,p}\in TERM(\tau)$. 
\smallskip
\noindent $\underline{\bf PROOF\, OF\, CLAIM\, ONE:}$\hfil \break

\noindent Let $\sigma\in TERM(\tau)$. We distinguish two cases: 

\leftline{{\bf Case One: }$\sigma_{p+1,p}=0$}
Here $\sigma-E_{p+1,p}$ is ineffective, so $$P(\sigma-E_{p+1,p})=0\,,$$
(and so the value chosen for its coefficient cannot affect the validity of 
(4.1.10)).

\leftline{{\bf Case Two: }$\sigma_{p+1,p}\ge1$}
Here $\sigma-E_{p+1,p}$ is effective. Now let us  list those terms in the 
three sums (4.1.9a,b,c) for which the assigned argument inside the 
symbol P, is precisely $\sigma-E_{p+1,p}$.

\noindent\underbar {Inside (4.1.9a):}\quad There occurs 
one such term, with coefficient given by (4.1.10a).

\noindent\underbar {Inside (4.1.9b):}\quad  Precisely one such 
term (i.e., containing $P(\sigma-E_{p+1,p})$) 
is assigned to each ordered pair
$(\sigma',k)$ with
$$\sigma'\in TERM(\tau),\quad i\le k<p$$
and such that
$$\sigma'-E_{p+1,k}-E_{p,k}=\sigma-E_{p+1,p}$$
---i.e., such that 
$$\sigma'=\sigma(k)'\buildrel \rm def \over =
\sigma-E_{p+1,p}+E_{p+1,k}-E_{p,k}\,.\eqno(4.1.11b)$$
To be explicit, to each such $(\sigma',k)$ corresponds
the term
$$T_b(\sigma',k)\buildrel \rm def \over =\langle 
\sigma' ;\tau \rangle (\sigma'_{p,k}+1)
P(\sigma-E_{p+1,p})\,.$$
Thus, if we set
$$F=\{k:i\le k<p\hbox{ and }\sigma(k)'\in TERM(\tau)$$
then the sum of the coefficients of $P(\sigma-E_{p+1,p})$
for all such $T_b(\sigma',k)$, is
$$\sum_F\{(\sigma(k)')_{p,k}+1\}\langle\sigma'(k);\tau\rangle
=\sum_F\sigma_{p,k} \langle\sigma(k)';\tau\rangle \,.$$
But there is no change in the value of this sum if we replace
$$\sum_F \hbox{ by }\sum_{k=i}^{p-1}$$
since by the conventions explained above, if $k$ is such that 
$$\sigma(k)'\not\in TERM(\tau)\,,$$ 
 then 
$$\langle\sigma(k)';\tau\rangle=0\,.$$
Thus the sum of the coefficients of the terms in question is precisely
the sum 
$$B_p(\sigma)=\sum_{k=i}^{p-1}\sigma_{p,k}\langle
\sigma(k)';\tau\rangle$$ 
given by eq.(4.1.10b).

\noindent\underbar {Inside (4.1.9c):}\quad  Essentially the same argument, 
 shows that   the terms in (4.1.9c) which involve
 $P(\sigma-E_{p+1,p})$, are precisely those of the form
$$T_c(\sigma'',k)\buildrel \rm def \over =-\langle
\sigma'' ;\tau \rangle (\sigma''_{k,p+1}+1)
P(\sigma-E_{p+1,p})\,,$$
with
$$\sigma''=\sigma(k)''\buildrel \rm def \over =
\sigma-E_{p+1,p}+E_{k,p}-E_{k,p+1}\,,\eqno(4.1.11c)$$
where $p+2\le k \le j$ and $$\sigma(k)''\in TERM(\tau)\,.$$
As before, these terms
have sum $C_p(\sigma)$ given by (4.1.10c).

\medskip
$\underline{\bf PROOF\,OF\,CLAIM\,TWO:}$

This is clear for the terms in (4.1.9a).

Consider next the terms in (4.1.9b).
Let
$$T=\langle \sigma';\tau \rangle
(\sigma'_{p,k}+1)P(\sigma'-E_{p+1,k}+E_{p,k})$$
be a non-zero term in (4.1.9b)---so, in particular, we have:
$$\sigma'\in TERM(\tau)\hbox{, and }\sigma'-E_{p+1,k}+E_{p,k}
\hbox{ is effective.}$$
Set
$$\sigma=\sigma'-E_{p+1,k}+E_{p,k}+E_{p+1,p}$$
Clearly $\sigma$ is effective, with $\sigma_{p+1,p}>0$. Since the weight-vector
$wt(\sigma)$ defined by eq.(1.10) is additive, and since 
$$\epsilon(E_{p+1,k}+E_{p,k}+E_{p+1,p})=0\,,$$
it follows that
$$\epsilon(\sigma)=\epsilon(\sigma')=r\alpha\,.$$
Hence $\sigma$ satisfies conditions II) and III) in Def.3.2.1.
It is also immediate from $k<p$ that $\sigma$ satisfies condition I)
in this definition.Hence $\sigma\in TERM(\tau)$,and so 
$$\sigma'=\sigma(k)',\hbox{ and the term T coincides with }
T_b(\sigma',k)\,,$$
as was to be proved.

The case of non-zero terms inside (4.1.9c) is precisely similar.

\centerline{\bf This completes the proof of Lemma 4.1.1.}
\bigskip
\par\noindent
The following simple lemma embodies an argument used in the
preceding, and
will find several further applications below.
\proclaim Lemma 4.1.2. Assume 
$$i\le p\le j,\, \sigma\in TERM(\tau),\, \sigma_{p+1,p}>0\,.$$
a)Suppose    also $ i\le k<p$ ; then we have:
$$\sigma(k)'\in TERM(\tau)\,\iff\sigma(k)'\hbox{ is effective,}$$
where $$\sigma(k)'\buildrel \rm def \over =
\sigma-E_{p+1,p}+E_{p+1,k}-E_{p,k}\,.$$ 
b)Suppose instead $p+2\le k<j $; then we have
$$\sigma(k)''\in TERM(\tau)\,\iff\sigma(k)''\hbox{ is effective,}$$
 where 
$$\sigma(k)''\buildrel \rm def \over =
\sigma-E_{p+1,p}+E_{k,p}-E_{k,p+1}\,.$$
\par\noindent
\smallskip
{\bf Proof:} If $i\le k<p$ then $-E_{p+1,p}+E_{p+1,k}-E_{p,k}$ has 
excess vector 0; since the excess vector $\epsilon$ is additive,
$$\epsilon(\sigma(k)')=\epsilon(\sigma)+0=d\alpha\,.$$
The assertion a) is now clear from the Remark in $\S3.2$. 
The proof of b) is similar.
\bigskip
\centerline{\bf Step Three\qquad   End-game: Proof $[E_{p,p+1},\gamma]=0$}
\bigskip
We are still assuming that $i\le p\le j-1$, and  that
$$\tau=(\lambda_i-\lambda_j,r,\sum_{k=1}^N l_k\lambda_k)$$
is a Verma triple---whence
$$l_i-l_j-i+j=r\,. \eqno(4.1.13)$$

Let $G_p$ denote the set of all  $\sigma$ in $ TERM(\tau)$ such that
$$\sigma-E_{p+1,p} \hbox{ is effective, i.e. such that }
\sigma_{p+1,p}\ge1 \eqno(4.1.14)$$
By Lemma 4.1.1, the desired result
$$[E_{p,p+1},\gamma]=0$$
will follow, if we show for every $\sigma$ in $G_p$, that 
the expression (4.1.10) equals 0.

Thus, to complete the proof of (4.1.3), it suffices to 
verify, for all $\sigma\in G_p$, that
$$A_p(\sigma)+B_p(\sigma)+C_p(\sigma)=0 \eqno(4.1.15)$$
where
$$\eqalign{
A_p(\sigma)&=\{l_p-l_{p+1}-col_p(\sigma)+col_{p+1}(\sigma)+1\}
\cdot \langle \sigma;\tau\rangle \cr
B_p(\sigma)&=\sum_{k=i}^{p-1}\sigma_{p,k}\langle
\sigma+E_{p+1,k}-E_{p,k}-E_{p+1,p};\tau \rangle \cr
&{\rm and}\cr
C_p({\sigma})&=-\sum_{k=p+2}^j \sigma_{k,p+1}\langle
\sigma+E_{k,p}-E_{k,p+1}-E_{p+1,p};\tau \rangle \cr
}$$ 

\bigskip
From here on, we shall write $R_k$ for
$$R_k(\sigma)=\sum_{l=i}^j\sigma_{k,l}=\sum_{l=i}^j\sigma_{l,k}
=col_k(\sigma)$$
(as given by Def.3.1.1 and by eqn. (4.1.6); and will also set
$$S_k=r-R_k\,.$$
\smallskip
The proof of (4.1.15) divides at this point into three cases, according as 
$$p=i,\; i<p<j-1 \hbox{  or } p=j-1\,.$$
\bigskip
\leftline{{\bf CASE I: }\qquad $\underline{p=i}$}  

Let $\sigma\in G_i$. Here (4.1.14) becomes $\sigma_{i+1,i}>0$
\ . Since $\sigma\in TERM(\tau)$,
$$col_i(\sigma)=r\, \hbox{\rm  and }\, col_{i+1}(\sigma)=
R_{i+1}\,,$$ 
so that
$$\eqalignno{
A_i(\sigma)&=(l_i-l_{i+1}-r+R_{i+1}+1)\cdot \langle
\sigma;\tau\rangle&(4.1.16)\cr
&=(l_i-l_{i+1}-S_{i+1}+1)\cdot r!\cdot\prod_{q=i+1}^{j-1}
R_q!S_q!{l_i-l_{i+1}-i+q \choose S_q} \cr}
$$
\bigskip
\par\noindent
$B_i(\sigma)=0$ in the present case, since the sum in (4.1.15b)
here becomes a sum over the empty indexing set $\{k:i\le k\le i-1\}$.
\bigskip
\noindent
We must next evaluate
$$C_i(\sigma)=-\sum_{k=i+2}^j \sigma_{k,i+1}\langle
\sigma(k)'';\tau\rangle$$
where, for
$$i+2\le k\le j \eqno(4.1.17)\,,$$
we set, as before,
$$\sigma(k)''=\sigma+E_{k,i}-E_{k,i+1}-E_{i+1,i}\,\,.$$
Now,
$$(\sigma(k)'')_{s,t}=\cases{
\sigma_{k,i}+1&if $(s,t)=(k,i)$\cr
\sigma_{k,i+1}-1&if $(s,t)=(k,i+1)$\cr
\sigma_{i+1,i}-1\ge 0&if $(s,t)=(i+1,i)$\cr
\sigma_{s,t}&otherwise\cr
}$$
and since $\sigma$ is effective, it follows that $\sigma(k)''$
is effective, if and only if $\sigma_{k,i+1}>0$.
It is now helpful (still assuming (4.1.17)), to divide into two sub-cases:
\hfil\break
\underbar{{\bf Sub-case} Ia:} \qquad $\underline{\sigma_{k,i+1}>0}$

Here, $\sigma(k)''$ is effective, hence, by Lemma 4.1.2, is in $ TERM(\tau)$.
Thus, Def.3.2.2 is applicable, and if we set
$$R''_q=R_q(\sigma(k)''),S''_q=r-R''_q \hbox{  for  }i<q<j\,,$$ 
we have
$$\langle\sigma(k)'';\tau\rangle=r!\cdot
\prod_{q=i+1}^{j-1}(R_q'')!(S_q'')!{l_i-l_q-i+q \choose S_q''}$$
We next compute the integers $R_q'',S_q''$ in this formula:

Since $\sigma(k)''\in TERM(\tau)$, 
$$R_q''=\sum_{s=1}^N (\sigma(k)'')_{q,s} 
=\sum_{s=i}^{q-1}(\sigma(k)'')_{q,s}$$
Thus,
$$R''_{i+1}=(\sigma(k)'')_{i+1,i}=\sigma_{i+1,i}+0-0-1=R_i-1,\, S''_{i+1}
=d-R''_{i+1}=S_i+1\,;$$
while for $i+2\le q\le j-1$,
$$R_q''=(\sigma(k)'')_{q+1,q}+\cdots+(\sigma(k)'')_{j,q}=
\sigma_{q+1,q}+\cdots+\sigma_{j,q}=R_q,\,S_q''=S_q\,.$$
Thus
$$\langle\sigma(k)'';\tau\rangle=(R_{i+1}-1)!(S_{i+1}+1)!
{l_i-l_{i+1}+1 \choose S_{i+1}+1}\cdot T \eqno(4.1.18)$$
where
$$T=r!\cdot \prod_{q=i+2}^{j-1} R_q!S_q!{l_i-l_q-i+q \choose S_q}\,.
\eqno(4.1.18a)$$
This of course trivially implies
$$\sigma_{k,i+1}\cdot\langle\sigma(k)'';\tau\rangle=
\sigma_{k,i+1}\cdot (R_{i+1}-1)!(S_{i+1}+1)!
{l_i-l_{i+1}+1 \choose S_{i+1}+1}\cdot T \eqno(4.1.19)$$
(This last deduction may seem less silly in a second.)

\underbar{{\bf Sub-case} Ib:} \qquad $\underline{\sigma_{k,i+1}=0}$
\ Here (4.1.18) is in general false, but (4.1.19) is still valid! Namely,
here $\sigma(k)''$ is {\bf not} effective, so by the conventions we are using,
$$\langle\sigma(k)'';\tau\rangle=0\,.$$
Thus, in general Def.3.2.2, i.e., (4.1.18), need not hold in this case---
but both sides of (4.1.19) are here 0.\hfil\break
(Note also that $\langle\sigma(k)'';\tau\rangle$ is independent of k 
in sub-case Ia), while it takes the (in general different) value
0 in sub-case Ib).)
\bigskip
Thus, we have proved that (4.1.19) holds for all k such that
(4.1.17) holds. Together with
$$R_{i+1}=\sum_{k=i+2}^j \sigma_{k,i+1}$$
we obtain
$$\eqalignno{
C_i(\sigma)&=-(\sum_{k=i+2}^j \sigma_{k,i+1}) \cdot (R_{i+1}-1)!
(S_{i+1}+1)!{l_i-l_{i+1}+1 \choose S_{i+1}+1}\cdot T \cr
&=-R_{i+1}!(S_{i+1}+1)!{l_i-l_{i+1}+1 \choose S_{i+1}+1}\cdot T \cr
}$$
Let us rewrite (4.1.16) as 
$$A_i(\sigma)=(l_i-l_{i+1}-S_{i+1}+1)\cdot R_{i+1}!
S_{i+1}!{l_i-l_{i+1}+1 \choose S_{i+1}}\cdot T$$
(with $T$ still given by (4.1.18a)) Thus we have
$$\displaylines{
A_i(\sigma)+B_i(\sigma)+C_i(\sigma)=  \cr
R_{i+1}!S_{i+1}!\left[(l_i-l_{i+1}-S_{i+1}+1){l_i-l_{i+1}+1 \choose
 S_{i+1}}-(S_{i+1}+1)!{l_i-l_{i+1}+1 \choose S_{i+1}}\right]T \cr  }$$
which is seen to vanish upon replacing, in the 
following elementary combinatorial
identity, $M$ by $l_i-l_{i+1}+1$ and  $N$ by $S_{i+1}$:
$$(M-N){M \choose  N}=(N+1){M \choose N+1} 
\eqno(4.1.20)\,,$$
(valid for M any complex number, N any non-negative integer.)
\bigskip
\centerline{\bf This proves (4.1.15) when $p=i$.}

\bigskip
\leftline{{\bf CASE II: }\qquad $\underline{i<p<j-1}$}
\smallskip
Let $\sigma\in G_p$. Thus,
$$\sigma\in TERM(\tau) \hbox{, and  }\sigma_{p+1,p}\ge 1 \,.$$ 
Both $p$ and $p+1$ are strictly between $i$ and $j$. We have
$$col_p(\sigma)=R_p \hbox{, and  }col_{p+1}(\sigma)=R_{p+1}\,,$$
and so 
$$A_p(\sigma)=(l_p-l_{p+1}-R_p+R_{p+1}+1)\cdot r! \cdot
\prod_{i<q<j}R_q!S_q!{l_i-l_q-i+q \choose S_q} \eqno(4.1.21)$$
\medskip
Consider next
$$B_p(\sigma)=\sum_{k=i}^{p-1} \sigma_{p,k}\langle \sigma(k)'
;\tau\rangle$$
where (as before)
$$\sigma(k)'=\sigma+E_{p+1,k}-E_{p,k}-E_{p+1,p}\,.$$
Its evaluation (which, with a few modifications, is quite similar to the 
evaluation just completed of $C_p(\sigma)$ in Case I) proceeds as follows:

Since $\sigma_{p+1,p}>0$, and $k<p$, and since $\sigma$ is effective, we
see that
$$\sigma(k)' \hbox{ is effective }\iff \;\sigma_{p,k}>0\,.$$
Accordingly, we now divide the study of $\sigma(k)'$
into two sub-cases, depending on whether
or not $\sigma_{p,k}$ is 0.
\bigskip
\underbar{{\bf Sub-case} IIa:} \qquad $\underline{\sigma_
{p,k}>0\;,i\le k<p}$

In sub-case IIa), $\sigma(k)'$ is effective, hence, by Lemma 4.2.2, is
in $ TERM(\tau)$. Thus, if we set
$$R_q'=R_q(\sigma(k)'),\,S_q'=S_q(\sigma(k)')=r-R_q'$$
then (in the present sub-case) 
$$ \langle \sigma(k)';\tau\rangle=r!\prod_{q=i+1}^{j-1}
(R_q')!(S_q')!{l_i-l_q-i+q \choose S_q'}$$
In this formula, the integers $R_q',S_q'$ associated with $\sigma(k)'$
coincide with the same integers for $\sigma$, except possibly for
$q=p,p+1$ or $k$; while
$$R_p'=R_p-1,R_{p+1}=R_{p+1},R_t'=R_t\hbox{; so  }
S_p'=S_p+1,S_{p+1}'=S_{p+1},S_t'=S_t\,.$$
Thus, the amplitude is here given by
$$ \langle \sigma(k)';\tau\rangle=(R_p-1)!(S_p+1)!{l_i-l_p-i+p
\choose S_p+1}\cdot T' \eqno(4.1.22)$$
where
$$T'=r!\cdot\prod_{\scriptstyle q\ne p \atop \scriptstyle i<q<j}
R_q!S_q!{l_i-l_q-i+q \choose S_q} \eqno(4.1.22a)$$
and so of course
$$\sigma_{p,k} \langle \sigma(k)';\tau\rangle=\sigma_{p,k}
(R_p-1)!(S_p+1)!{l_i-l_p-i+p \choose S_p+1}\cdot T' \eqno(4.1.23)$$
\medskip
\underbar{{\bf Sub-case} IIb:} \qquad $\underline{\sigma_
{p,k}=0\;,i\le k<p}$
\medskip
(4.1.23) remins valid in the present sub-case, since both sides are 0.
\bigskip
Thus, in Sub-cases IIa) and IIb) alike, (4.1.23) is valid, and since
$$R_p=\sum_{k=i}^{p-1}\sigma_{p,k}$$
it follows that
$$\eqalignno{
B_p(\sigma)&=\sum_{k=i}^{p-1}\sigma_{p,k}(R_p-1)!(S_p+1)!
{l_i-l_p-i+p \choose S_p+1}T' &(4.1.24) \cr
&=R_p!(S_p+1)!{l_i-l_p-i+p \choose S_p+1}T'\cr
}$$
\bigskip
We next turn to the computation of
$$C_p(\sigma)=-\sum_{k=p+2}^j \sigma_{k,p+1}
\langle \sigma(k)'';\tau\rangle$$
where now 
$$\sigma(k)''=\sigma+E_{k,p}-E_{k,p+1}-E_{p+1,p}\,.$$
This computation is similar to the preceding one of $B_p(\sigma)$,
with the following modifications:

Assume $p+2\le k\le j,\; \sigma\in G_p$.Then
$$\sigma(k)'' \hbox{ is effective }\iff \sigma_{k,p+1}>0 \,,$$
and we have two sub-cases to consider, according as $\sigma_{k,p+1}>0$
 or not.

If $\sigma_{k,p+1}>0$ then (by Lemma 4.1.2) $\sigma(k)''$ is in
$ TERM(\tau)$, so we may use Def.3.2.2 to compute the amplitude
$\langle \sigma(k)'';\tau\rangle$. Since $i<s<j$, one readily verifies that
$$R_q(\sigma(k)'')=\cases{R_{p+1}-1& if $q=p+1$ \cr
R_q& if $q\ne p+1$ \cr}\,;\,S_q(\sigma(k)'')=\cases{
S_{p+1}+1 & if $q=p+1$ \cr S_q& if $q\ne p+1$ \cr}$$
whence we obtain
$$\langle \sigma(k)'';\tau\rangle=(R_{p+1}-1)!(S_{p+1}+1)!
{l_i-l_{p+1}-i+p+1 \choose S_{p+1}+1}\cdot T'' \eqno(4.1.25)$$
where
$$T''=r!\prod_{\scriptstyle q\ne p+1 \atop \scriptstyle i<q<j}
R_q!S_q!{l_i-l_q-i+q \choose  S_q} \eqno(4.1.25a) $$

Thus, whether $\sigma_{k,p+1}>0$ or not, 
$$\sigma_{k,p+1}\cdot \langle \sigma(k)'';\tau\rangle=
\sigma_{k,p+1}\cdot (R_{p+1}-1)!(S_{p+1}+1)!
{l_i-l_{p+1}-i+p+1 \choose S_{p+1}+1}\cdot T'' \eqno(4.1.26)$$ 
is always valid. Together with 
$$R_{p+1}=\sum_{k=p+2}^j \sigma_{k,p+1}\,,$$
this yields (as in the preceding discussion of $B_p(\sigma)$)
$$C_p(\sigma)=-R_{p+1}!(S_{p+1}+1)!{l_i-l_{p+1}-i+p+1 \choose 
 S_{p+1}+1}\cdot T'' \eqno(4.1.27) $$

It is now  convenient to introduce the following common divisor of the
right-hand sides of eqns.(4.1.21),(4.1.24) and (4.1.27):
$$T=r!\cdot R_p!S_p!(R_{p+1})!(S_{p+1})!\cdot
\prod_{\scriptstyle q\ne p,q\ne p+1 \atop \scriptstyle
i<q<j}R_q!S_q!{l_i-l_q-i+q \choose  S_q} \eqno(4.1.28) $$
Then we have
$$A_p(\sigma)+B_p(\sigma)+C_p(\sigma)=TU \eqno(4.1.29)$$
with
$$\eqalignno{U&=(l_p-l_{p+1}-R_p+R_{p+1}+1){l_i-l_p-i+p
\choose S_p} &(4.1.29a) \cr
&+(S_p+1){l_i-l_p-i+p \choose S_p+1}{l_i-l_{p+1}-i+p+1
\choose S_{p+1}} \cr
&-(S_{p+1}+1){l_i-l_p-i+p \choose S_p}{l_i-l_{p+1}-i+p+1
\choose S_{p+1}+1} \cr }$$
Thus, to prove (4.1.15), it suffices in the present case to verify that
$U=0$.

For this purpose, first observe 
$$l_p-l_{p+1}-R_p+R_{p+1}+1=-(l_i-l_p-i+p-S_p)+
(l_i-l_{p+1}-i+p+1-S_{p+1})$$
so that we may write
$$U={l_i-l_p-i+p \choose S_p}U_1+{l_i-l_{p+1}-i+p+1 \choose
 S_{p+1}}U_2$$
with 
$$U_1=(l_i-l_{p+1}-i+p+1-S_{p+1})
{l_i-l_{p+1}-i+p+1 \choose S_{p+1}}-(S_{p+1}+1){l_i-l_{p+1}-i+p+1 
\choose S_{p+1}+1} $$
and
$$U_2=-(l_i-l_p-i+p-S_p){l_i-l_P-i+p \choose S_p}
+(S_p+1){l_i-l_p-i+p \choose S_p+1}$$
and then observe that $U_1$ and $U_2$ both vanish because of (4.1.20).
\bigskip
\centerline{\bf This proves (4.1.15) when $i<p<j-1$.}
\bigskip
\leftline{{\bf CASE III: }\qquad $\underline{p=j-1}$}
\smallskip
{\bf Remark:} Up to this point, the argument has made no use 
of the Verma condition. In other words, except for the single value 
$j-1$ for $p$, (4.1.3)
holds with no restrictions on the value of $\lambda$. Thus, it can be 
predicted that we shall need to utilize the Verma condition (4.1), 
in the present (final remaining) case,
i.e. to utilize
$$l_i-l_j-i+j=r\,.$$ 
\bigskip
Let us begin by noting that
$$col_p(\sigma)=col_{j-1}(\sigma)=R_{j-1},\,col_{p+1}(\sigma)=
col_j(\sigma)=0\,,$$
so
$$A_{j-1}(\sigma)=(l_{j-1}-l_j-R_{j-1}+1)\cdot r!\cdot
\prod_{i<q<j}R_q!S_q!{l_i-l_q-i+q \choose S_q} \eqno(4.1.30)$$
\medskip
Next, we note that in the present case,
$$C_p(\sigma)=0\,,$$
since the sum $\sum_{k=p+2}^j$ on the right side of (4.1.15c), is
here extended over the empty set $\{k:j+1\le k\le j\}$.
\medskip
Finally, we must evaluate
$$B_p(\sigma)= B_{j-1}(\sigma)=\sum_{k=i}^{j-2}\sigma_{j-1,k}
\langle \sigma(k)';\tau\rangle$$
with $\sigma(k)'$ defined (for the present value of p, and for
$ i\le k\le j-2$) by
$$\sigma(k)'=\sigma-E_{j,j-1}+E_{j,k}-E_{j-1,k}
\hbox{  (for $ i\le k\le j-2$)}$$

\noindent
{\bf CLAIM:}For $ i\le k\le j-2$,
$$\sigma_{j-1,k}\langle \sigma(k)';\tau\rangle=
\sigma_{j-1,k}\cdot (S_{j-1}+1){l_i-l_{j-1}-i+j-1 \choose
 S_{j-1}+1}\cdot T \eqno(4.1.31) $$
where
$$T=r!R_{j-1}!S_{j-1}!\prod_{q=i+1}^{j-2}R_q!S_q!{l_i-l_q-i+q \choose S_q}
\eqno(4.1.31a)$$

The proof is much the same as that used before, to derive the value
obtained for $C_i(\sigma)$ in Case I, and used in Case II to derive
$B_p(\sigma)$ (eqn. 4.1.24), and $C_p(\sigma)$ (eqn.4.1.27). Namely:

If $\sigma_{p-1,k}=0$, then both sides of (4.1.31) are 0. Otherwise,
 $\sigma(k)'$ is effective, hence lies in $ TERM(\tau)$ (by a last
use of Lemma 4.2.2), so we may use Def.3.2.2 to compute the amplitude
$\langle \sigma(k)';\tau\rangle$. Noting for this purpose, that if
$i<s<j$, then we have
$$R_s(\sigma(k)')=\cases{R_{j-1}-1&if $s=j-1$ \cr
R_s&if $s<j_1$ \cr }\hbox{ ;so }S_s(\sigma(k)')=\cases{
S{j-1}+1&if $s=j-1$ \cr S_s & if $s<j-1$ \cr }$$
it follows immediately that
$$\langle \sigma(k)';\tau\rangle=(S_{j-1}+1){l_i-l_{j-1}-i+j-1 
\choose S_{j-1}+1}\cdot T$$
\centerline{\bf This completes the proof of Claim (4.1.31).}

Combining (4.1.31) with 
$$R_{j-1}=\sum_{k=i}^{j-2}\sigma_{j-1,k}$$
now yields the desired evaluation
$$B_{j-1}(\sigma)=(S_{j-1}+1){l_i-l_{j-1}-i+j-1 \choose
S_{j-1}+1}\cdot T$$
Let us also rewrite (4.1.30) as
$$A_{j-1}(\sigma)=(l_{j-1}-l_j-R_{j-1}+1)
{l_i-l_{j-1}-i+j-1 \choose S_{j-1}}\cdot T $$

Then we get
$$\displaylines{
A_{j-1}(\sigma)+B_{j-1}(\sigma)+C_{j-1}(\sigma)= \cr
\left[(l_{j-1}-l_j-R_{j-1}+1){l_i-l_{j-1}-i+j-1
\choose S_{j-1}}+(S_{j-1}+1){l_i-l_{j-1}-i+j-1 \choose
S_{j-1}+1)} \right]\cdot T \cr
}$$
Thus, in order to prove (4.1.15), in the present case, it is sufficient 
to verify that the following expression U vanishes:
$$U=(l_{j-1}-l_j-R_{j-1}+1){l_i-l_{j-1}-i+j-1
\choose S_{j-1}}+(S_{j-1}+1){l_i-l_{j-1}-i+j-1 \choose
S_{j-1}+1)}$$
Using (4.1.20), we obtain 
$$\eqalign{
U&={l_i-l_{j-1}-i+j-1 \choose S_{j-1}}[(l_{j-1}-l_j-R_{j-1}+1)+
(l_i-l_{j-1}-i+j-1-S_{j-1})] \cr
&={l_i-l_{j-1}-i+j-1 \choose S_{j-1}}[l_i-l_j-i+j-(R_{j-1}+S_{j-1})] \cr
&={l_i-l_{j-1}-i+j-1 \choose S_{j-1}}[l_i-l_j-i+j-r] \cr
}$$
which indeed vanishes whenever the Verma condition (4.1)---i.e.
$$l_i-l_j-i+j=r$$
---is satisfied. 

We have thus proved eqn.(4.1.15), hence also eqn.(4.1.3), in all three
cases.
\bigskip
\centerline {\bf This completes the proof of the assertion A1. }
\bigskip
\centerline{\bf \S4.2 Proof that $\gamma\cdot v_{\lambda}$ 
has Weight $\lambda - r\alpha$}
\medskip
We next turn to the proof of A2, which asserts that the element
 $\gamma\cdot v_{\lambda}$ in the ${\frak sl}_N$-module
${\cal V}_{\lambda}$ has weight $\lambda - r\alpha$.

Now, for every $N$-shift $\sigma$, it follows from eqn.(1.11) in
\S1.3, that $P(\sigma)v_{\lambda}$ has weight
$$wt(\sigma)+wt(v_\lambda)=wt(\sigma)+\lambda \,,$$
where the weight $ wt(\sigma )$ is given by eqn.(1.10).
Since $\gamma$ is, by definition, a ${\Bbb Z}$-linear combination of
$$ \{P(\sigma):\sigma \in TERM(i,j,r)\}$$
the desired assertion A2  will be an immediate consequence of the following
lemma:
\smallskip
\proclaim Lemma 4.2.1. Let $\sigma\in TERM(i,j,r)$;then
$$wt(\sigma)=-r\cdot (\lambda_i-\lambda_j)$$
\par 

{\bf PROOF:}
By hypothesis, $\sigma$ satisfies the three conditions I),II),III)
of Def.3.2.1. It must then be proved that, for all k between 1 and N 
inclusive,
$$wt_k(\sigma)=-r\cdot(\lambda_i-\lambda_j)(E_{k,k})$$
i.e., that
$$\sum_{l=1}^N \sigma_{k,l}-\sum_{l=1}^N \sigma_{l,k}=r\cdot(\delta_{j,k}
-\delta_{i,k}) \eqno(4.2.1)$$
There are four cases to check:\par

\noindent {\bf CASE 1: } $\underline {k=i}$ \hfil \break
Here the first sum in (4.2.1) vanishes, by condition I) of
Def.3.2.1, so we are left with  
$$\sum_{l=1}^N \sigma_{i,l}= r$$
which is precisely condition II).\par

\noindent {\bf CASE 2: } $\underline {i<k<j}$ \hfil \break 
Here the equation to be proved is
$$\sum_{l=1}^N \sigma_{k,l}-\sum_{l=1}^N \sigma_{l,k}=0$$
which holds (for all k in the given range) by condition III).\par

\noindent {\bf CASE 3: } $\underline {k=j}$ \hfil \break 
Here the second sum in (4.2.1) vanishes by Condition I), and
we are left to prove:
$$\sum_{l=i}^{j-1}\sigma_{j,l}=r \eqno(4.2.2) $$
For $i \le s <j$ let us set
$$C(s)=\sigma_{s+1,s}+\sigma_{s+2,s}+\cdots+\sigma_{j,s}
+\sum_{i\le p<s<q\le j}\sigma_{q,p}  $$
Note that 
$$C(i)=\sum_{i<k\le j}\sigma_{k,i}=r$$
by condition I), while for $i< s<j$,
$$C(s)-C(s-1)=\sum_{i \le k <s}\sigma_{s,k}-
\sum_{s<k \le j}\sigma_{s,k}$$
which equals 0 by Condition III) of Def.3.2.1.Hence $C(j-1)=r$, which
is the same as (4.2.2).\par

\noindent {\bf CASE 4:}$\underline {k<i \hbox{ or } k>j}$ \hfil \break 
Here both of the sums occurring in (4.2.1) are empty.

{\bf This completes the proof of (4.2.1), and hence of Assertion A2.}

\bigskip
\centerline{\bf \S4.3  Proof that $\gamma$ satisfies the Shapovalov 
normalization condition VS2)}

\medskip
Let us define an $N$-shift $\sigma$ to be {\it lower-triangular}
if it has this property:
$$ \hbox{for all }k,l \in \underline N\,,\;\sigma_{k,l}\ne 0 \Rightarrow
k>l \,.$$  
Denote by $(\Pi^N)^{LT}$ the set of all such.
(Note that all the $P(\sigma)$ involved in the formula (4.1) for
$\gamma$ have lower-triangular $\sigma$.)

Choose (arbitrarily) a total ordering $<<$ for $\Delta^+$, say
$$\alpha_1 << \cdots << \alpha_m$$
where
$$\alpha_s=(i(s),j(s)) \hbox{ with } i(s) < j(s)\, (1 \le s\le m)$$ 

This in turn induces a total ordering (which by a slight abuse of notation
will also be denoted by $<<$ ) on the set of all lower-triangular
$N$-shifts, defined as follows:

To each lower-triangular  $N$-shift $\sigma$ assign the ordered $m$-tuple
$$\langle \sigma \rangle \buildrel \rm def \over{=}
(\sigma_{j(1),i(1)},\sigma_{j(2),i(2)},\cdots,\sigma_{j(m),i(m)}) $$
Given a second lower-triangular  $N$-shift $\tau$, we shall say that $\sigma$ 
precedes $\tau$ in the given total ordering,
$$\sigma << \tau\,,$$
if and only if $\langle \sigma \rangle $ precedes
$\langle \tau \rangle $ in the usual lexicographic ordering on
ordered $m$-tuples of non-negative integers.

The {\it weight} $W(\sigma)$ of any $N$-shift $\sigma$ is defined
to be
$$W(\sigma)\buildrel \rm def \over{=} \sum_{k,l\in \underline N}
\sigma_{k,l}\,.$$    

As in \S3.1, we assign to every lower-triangular
$N$-shift $\sigma$, the basis-vector
$$F_{\sigma}:=(E_{j(1),i(1)})^{\sigma(j(1),i(1))}\cdot\cdots
\cdot (E_{j(m),i(m)})^{\sigma(j(m),i(m))}$$  
in the PBW-basis 
$${\frak B}(<<)=\{F_{\sigma} \bigl | \sigma \in (\Pi^N)^{LT}\} \eqno(4.3.1) $$
for ${\frak A}(N_-)$ determined by $<<$.

We take the usual filtration for the enveloping algebra
 ${\frak A}(N_-)$, whereby the $s$-th filtration-level ${\frak A}(N_-)_s$ is
the ${\Bbb C}$-span of all products of $\le s$ elements $E_{i,j}$ with
$i>j$. Note that then ${\frak A}(N_-)_s$ has basis consisting of
$$\{F_{\sigma} \bigl | \sigma \hbox{ is lower-triangular
of weight} \le s \}$$ 

\medskip
\proclaim Lemma 4.3.1. Let $\sigma$ be a lower-triangular $N$-shift.    
Let $1\le k<l\le N$.Then we may write
$$E_{l,k}F_{\sigma}=F_{E_{l,k}+\sigma}+
\sum_{p=1}^L c_p \cdot F_{\sigma_p}$$
where $L$ is a non-negative integer (which may be zero), each
$c_p$ is a complex number, and 
$$\sigma_1,\cdots,\sigma_L$$
are lower-triangular $N$-shifts of weight smaller than that of $\sigma$.

\leftline {{\bf PROOF:}} 
\noindent Immediate.
\medskip

\proclaim Lemma 4.3.2. Let $\sigma$ be a lower-triangular $N$-shift.
Express the Weyl polarization $P(\sigma)$ as a ${\Bbb C}$-linear
combination of the basis ${\frak B}(<<)$ for ${\frak A}(N_-)$, say:
$$ P(\sigma)=\sum \{C_{\tau} \cdot F_{\tau} \bigl | \tau \in 
(\Pi^N)^{LT}\} \eqno(4.3.2) $$
(where all $C_{\tau}$ are complex numbers).
Then 
$$C_{\sigma}={1 \over (\sigma)!}\,,$$
and
$$C_{\tau}\ne 0 \Rightarrow W(\tau)<W(\sigma)\,.$$

\leftline{{\bf PROOF:}} We argue by induction on $W=W(\sigma)$:

If W=1, then for suitable $k,l$ in $\underline N$, we have
$$k<l, \sigma=E_{l,k}\,.$$
Hence $P(\sigma)=D_{l,k}=F_{E_{l,k}}$, ${\sigma}!=1$, and the assertion is 
clear in this case.

Next, assume that $W>1$, and that the Lemma to be proved holds
for all lower-triangular $N$-shifts of weight $< W$. 

It follows immediately from the induction hypothesis, that
if $\tau$ is a lower-triangular $N$-shift of weight $W'<W$, then
$P(\tau)$ lies in  ${\frak A}(N_-)_{W'}$.

There exist integers $k,l$ with
$$1 \le k<l\le N,\,\sigma_{l,k}> 0 \,.$$
By Cor.2.3.2, there then exist lower-triangular N-shifts 
$$\sigma_1,\cdots,\sigma_L$$
(with $L \ge 0$), each of weight $W-1$,
 and positive integers $m_1,\cdots,m_L$, such that
$$\sigma_{l,k}\cdot P(\sigma)=E_{l,k}P(\sigma -E_{i,j})-
\sum_{s=1}^L m_s P(\sigma_s)\,. \eqno(4.3.3)$$ 
By the induction hypothesis applied to the $N$-shift
$\sigma-E_{l,k}$ of weight $W-1$, we may write
$$P(\sigma-E_{l,k})={1 \over (\sigma-E_{l,k})!}F_{\sigma-E_{l,k}}
+(  \hbox { an element in } {\frak A}(N_-)_{W-2}) \,.$$
Hence, using Lemma 4.3.1, it follows that 
$$E_{l,k}P(\sigma-E_{l,k})={1 \over (\sigma-E_{l,k})!}F_{\sigma}
+(  \hbox { an element in } {\frak A}(N_-)_{W-1}) \,.$$
Since 
$$(\sigma-E_{l,k})!={1 \over \sigma_{l,k}}\cdot (\sigma)! \,,$$
we may rewrite (4.3.3) as
$$P(\sigma)={1 \over {\sigma}!}F_{\sigma}+
(  \hbox { an element in } {\frak A}(N_-)_{W-1})$$
which completes the induction, and so the proof of Lemma 4.3.2. 
\bigskip

Our present purpose is the proof of assertion A3, i.e., the proof that 
$\gamma$ satisfies the Shapovalov normalization property
VS2) explained in \S3.1. Let $\sigma_0$ denote the $N$-shift
$$\sigma_0=r(E_{i+1,i}+\cdots+E_{j-1,j})\,. \eqno(4.3.4)$$
The assertion to be proved is, that when
$\gamma$ is expressed as a ${\Bbb C}$-linear combination of the basis
(4.3.1) for ${\frak A}(N_-)$, the coefficient of the basis vector
$F_{\sigma_0}$ is precisely 1.
\smallskip
To see this, let us recall eqn.(4.1), which defined $\gamma$ as:
$$\gamma\buildrel\rm def\over =\sum_{\sigma \in {\cal T}(\tau)}
\langle \sigma;\tau \rangle P(\sigma) \eqno(4.3.5)$$ 
where 
$${\cal T}(\tau)=TERM(i,j,r)\,,$$
is the set of lower-triangular $N$-shifts 
subordinate to $(\lambda_i-\lambda_j,r)$, as
furnished by Def.3.1.1. In particular, the three conditions of
Def.3.1.1 are satisfied by the $N$-shift (4.3.4), i.e. ${\cal T}(\tau)$
contains $\sigma_0$.

{\bf CLAIM:} $\sigma_0$ has strictly larger weight than any other member 
of ${\cal T}(\tau)$.

{\bf PROOF:}If $\sigma_1$ is any element of ${\cal T}(\tau)$, and  $\sigma_1$
is distinct from $\sigma_0$, then we may write
$$\sigma_1=E_{l,k}+\sigma_1'$$
with 
$$i\le k<l\le j,\,l-k>1,\;\hbox{ and $\sigma_1'$ effective. }$$
Then the N-shift
$$E_{l,k+1}+E_{k+1,k}+\sigma_1'$$
lies in ${\cal T}(\tau)$, and has weight larger by 1 than that of
$\sigma_1$---which completes the proof of the Claim.

Combining the result just proved with Lemma 4.3.2 and eqn.(4.3.5), we
obtain:
$$\gamma={{\langle \sigma_0;\tau \rangle} \over {\sigma_0}!} F_{\sigma_0}+
(  \hbox { an element in } {\frak A}(N_-)_{r(i-j)-1})$$     
i.e., there exist distinct lower-triangular 
$$\sigma_{1},\sigma_{2},\cdots,\sigma_{L}\,,$$
all of weight strictly less than that of $\sigma_0$, 
 and complex numbers $C_s$, such that
$$\gamma={{\langle \sigma_0;\tau \rangle \over {\sigma_0}!}} F_{\sigma_0}+ 
\sum_{s=1}^L C_s F_{\sigma_s}\eqno(4.3.6)$$
Let us now appeal to Def.3.2.2 in order to 
compute $\langle \sigma_0;\tau \rangle$.
Here we must replace each $K_\sigma(k)$ in eqn.(3.2.2) by $r$, thus 
obtaining
$$\langle \sigma_0;\tau \rangle=(r!)\cdot \prod_{k=i+1}^{j-1}
[r!0!{l_i-l_k \choose 0}]= (r!)^{j-i} $$
But clearly eqn.(4.3.4) implies that
$$(\sigma_0)!=(r!)^{j-1}$$
Hence, we see that the coefficient of $F_{\sigma_0}$ in the right-hand
side of eqn.(4.3.6) is 1, as was to be proved.

This completes the proof of Assertion A3, hence of
Theorem 3.3.4, and so also of all the assertions in \S3.

\def\bib#1{\itemitem{[#1]}}
\hyphenation{U-ni-ver-si-tet Ma-te-ma-ti-ske}

\noindent
\centerline{References}
\bigskip
{\narrower

\bib{Akin1} K.AKIN, On Complexes Relating the Jacobi-Trudi Identity
with the Bernstein-Gelfand-Gelfand Resolution I. J.\ of Alg.\ {\bf 117}
(1988), 494--503

\bib{Akin2} K.AKIN, On Complexes Relating the Jacobi-Trudi Identity
with the Bernstein-Gelfand-Gelfand Resolution II. J.\ of Alg.\ {\bf 152}
(1992), 417--426 

\bib{BGG} I.N.BERNSTEIN, I.M.GEL'FAND and S.I.GEL'FAND, Differential
operators on the base affine space and a study of ${\frak g}$-modules,
in `Lie Groups and their Representations', Proc. of the summer school on
group representations, Bolyai Ja\'nos Math.\ Society, Budapest, (1971)
, Wiley (New York--Toronto 1975) 21--64 

\bib{Cap1}A.CAPELLI, \"Uber die Zur\"uckf\"urung der Cayley'schen
Operation $\Omega$ auf gew\"ohnliche Polar-Operationen, Math.Ann.
{\bf 29}(1887), 331--338

\bib{Cap2}A.CAPELLI, Sur les op\'erations dans la th\'eorie des formes 
alg\'ebriques, Math.Ann.{\bf 37}(1890),1--37

\bib{Cap3}A.CAPELLI, Lezioni sulla Theoria delle Forme Al\-ge\-bri\-che,
\hfil\break  Pel\-ler\-ano,Na\-po\-li,1902

\bib{CL}R.W.CARTER AND G.LUSZTIG, On the modular representations of the
general linear and symmetric groups. Math.Z.{\bf 136}(1974)193--242

\bib{Doty} S.R.DOTY, The Symmetric Algebra and Representations of 
General Linear Groups,  Proceedings of the Hyderabad Conference on 
Algebraic Groups, (Dec.1989), 123--150 

\bib{Doty2} S.R.DOTY, Resolutions of $B$ Modules, Indag.Mathem.,N.S.,
{\bf 5(3)},267--283 (Sept. 26, 1994)

\bib{Fra} J.FRANKLIN, Homomorphisms between Verma Modules in 
Characteristic p. J.\ of Alg.\ {\bf 112}(1988)58--85

\bib{FH} W.FULTON and J.HARRIS, Representation Theory, A First Course,
{\bf 129}, Graduate Texts in Mathematics, Springer Verlag, 1991 

\bib{Greene}GREENE

\bib{Las}A.LASCOUX, Syzygies des vari\'et\'es d\'eterminentales,
Advances in Math. {\bf301}, (1978), 202--237

\bib{Mal}M.MALIAKAS,Resolutions and parabolic Schur algebras.
 J.\ Alg.\ {\bf 180},\hfil\break (1996),679---690

\bib{Mathas} A.MATHAS, Iwahori-Hecke Algebras and Schur Algebras of the 
Symmetric Group, AMS University Lecure Series {\bf 15}(1991)

\bib{MFF}F.G.MALIKOV;B.L.FEIGIN;D.FUKS. 
 Singular vectors in Verma
modules over Kac - Moody algebras. (Russian) Funkts. Anal. i
Prilozhen. {\bf 20} (1986), no. 2, 25--37, 96.

\bib{Nielsen} A.NIELSEN, Tensor Functors of Complexes, Aarhus 
Universitet, Matematiske Institut, Preprint Series 77/78, No.15

\bib{Shap} N.N.SHAPOVALOV, On a bilinear form on the universal
enveloping algebra of a complex semisimple Lie algebra, Functional
 Anal.\ Appl.\ .{\bf 6}, (1972), 307--312

\bib{Umeda} T.UMEDA, The Capelli Identities, A Century 
After, Amer.\ Math.\ Soc.\ Transl.
 {\bf 2 Vol.183}, 1998

\bib{Verma1} D.-N. VERMA, Structure of Certain Induced Representations
of Complex Semisimple Lie Algebras, Yale University Doctoral Dissertation,
 (1966)

\bib{Verma2}  D.-N. VERMA, Structure of certain induced representations
of complex semisimple Lie algebras, Bull.A.M.S. {\bf 74}(1968), 160--166

\bib{Weyl}H.WEYL,  The Classical Groups, Their Invariants and
Representations, Princeton University Press, 1946

\bib{Wood}D.J.WOODCOCK, Comm.\ Alg.\ {\bf 22, no.5},(1994),1703--1721

\bib{Young} A.YOUNG, On Quantitative Substitutional Analysis I, PLMS
{\bf 33} (1900) 97--146

\bib{Zel} A.ZELEVINSKY, Resolvents, Dual pairs, and Character Formulas,
\hfil\break
 Functional  Anal.\  Appl.\ , {\bf 21} (1987), 152--154
\smallskip}

\bye